\documentclass[a4paper, reqno, 12pt]{amsart}
\usepackage{tikz-cd}
\usepackage[all]{xy}
\usepackage{amssymb}
\usepackage{amsmath,amscd,color}
\usepackage[all]{xy}
\usepackage{mathrsfs}
\usepackage{tensor}
\usepackage{stackrel}
\usepackage{hyperref}
\usepackage{mathrsfs}
\usepackage{hyperref}
\usepackage{mathtools}
\usepackage[margin=1.2 in]{geometry}
\usepackage{leftindex}
\usepackage{enumitem}

\definecolor{tocolor}{rgb}{.1,.1,.5}
\definecolor{urlcolor}{rgb}{.2,.2,.6}
\definecolor{linkcolor}{rgb}{.1,.1,.6}
\definecolor{citecolor}{rgb}{.6,.2,.1}
\hypersetup{ colorlinks=true, urlcolor=urlcolor,
	linkcolor=linkcolor, citecolor=citecolor}
\definecolor{darkgreen}{rgb}{0.0, 0.5, 0.0}

\providecommand{\U}[1]{\protect\rule{.1in}{.1in}}
\newtheorem{theorem}{Theorem}[section]
\newtheorem*{theorem*}{Theorem}
\newtheorem*{claim*}{Claim}

\newtheorem{corollary}[theorem]{Corollary}

\newtheorem{proposition}[theorem]{Proposition}

\newtheorem*{proposition*}{Proposition}
\newtheorem{lemma}[theorem]{Lemma}

\newtheorem*{question}{Question}

\theoremstyle{definition}
\newtheorem{definition}[theorem]{Definition}

\newtheorem{example}[theorem]{Example}
\newtheorem{remark}[theorem]{Remark}
\numberwithin{equation}{section}

\newcommand{\mf}[1]{\mathfrak{#1}}

\newcommand{\del}{\partial}
\newcommand{\toto}{\rightrightarrows}

\newcommand{\RR}{\mathbb{R}}                    

\newcommand{\LL}{\mathcal{L}}                   
\newcommand{\GG}{\mathcal{G}}                   
\renewcommand{\:}{\colon}                       
\renewcommand{\d}{{\mathrm{d}}}                   
\newcommand{\D}{{\mathrm{D}}}                   
\newcommand{\DD}{\mathcal{D}}                   
  
\newcommand{\Cc}{\mathcal{C}}                    

\newcommand{\Hom}{\operatorname{Hom}}       
\newcommand{\id}{\operatorname{id}}         
\newcommand{\im}{\operatorname{im}}         
\newcommand{\Der}{\operatorname{Der}}       
\newcommand{\UD}{\Breve{\mathrm{D}}}                        

\DeclareMathOperator{\Diff}{Diff}           
\newcommand{\loc}{\mathrm{loc}}             
\newcommand{\Diffloc}{\Diff_{\loc}}         
\newcommand{\dom}{\operatorname{dom}}       
\newcommand{\pullback}[2]{{}_{#1}\kern-\scriptspace{\times}_{#2}}                   

\newcommand{\onabla}{{\overline{\nabla}}}   

\begin{document}
	\title{Notes on relative algebroids}

	\author{Wilmer Smilde}

	\address{Department of Mathematics, University of Illinois at Urbana-Champaign, 1409 W. Green Street, Urbana, IL 61801 USA}
	\email{wsmilde2@illinois.edu}
	
	\date{\today}
\begin{abstract}
Relative algebroids provide a framework that unifies Lie algebroids with partial differential equations. In this set of notes, we explain how relative algebroids arise from geometric problems, and give an introduction to their structural theory. We also discuss their relation to and relevance for partial differential equations with symmetry. 
\end{abstract}

\thanks{This work was partially supported by NSF grant DMS-2303586 and UIUC Campus Research Board Award RB25014.}

\maketitle
\setcounter{tocdepth}{2}
\tableofcontents
\tableofcontents
\section*{Introduction}
These notes provide an introduction to the theory of relative algebroids. They are meant complement the main reference \cite{FernandesSmilde2025} as a more readable and digestible exposition the framework. They cover the main ideas and core constructions while minimizing technical language. The focus of this article is restricted to algebroids relative submersions, as opposed to more general algebroids relative to foliations discussed in \cite{FernandesSmilde2025}. Central topics of the theory presented here are prolongation and formal integrability of relative algebroids, as well as their relation with PDEs and symmetries. There is also more emphasis on examples compared to \cite{FernandesSmilde2025}, and some (simplified versions of) results in \cite{Smilde2025} are announced. The trade-off is that several technical results are skipped. Most notably, the fundamental but technical notion of a tableaux of derivations is completely circumvented. Nevertheless, we hope that this set of notes makes the ideas surrounding the theory of relative algebroids accessible.

The notes are set up as follows. In Section \ref{sec:MotivatingExamplesAndOverview} we discuss several examples of the appearance of algebroids (relative or not) in existence and classification problems for geometric structures. We also give an overview of the role of Lie algebroids and relative algebroids in understanding such problems.

In Section \ref{sec:RelativeDerivations}, we dive into the structural theory of relative derivations. This builds the necessary technical language for the later sections.

In Section \ref{sec:RelativeAlgebroids}, we introduce the concept of a relative algebroid, and discuss through computational examples the notions of prolongation and formal integrability. At the end, we revisit an example from Section \ref{sec:MotivatingExamplesAndOverview} to make the theory concrete.

In Section \ref{sec:RelativeAlgebroidsAndPDEs} we make clear the connection between partial differential equations (PDEs) and relative algebroids. This is important because it explains how the known theory of formal PDEs fits inside the theory of relative algebroids. We also discuss invariance of the relative algebroid under symmetries, and give an application of taking a quotient of a PDE by symmetries, resulting in a relative algebroid assocatied to relative distributions.

\subsubsection*{Acknowledgments}
These lecture notes grew out of several talks and mini-courses on relative algebroids, most notably a talk at the \textit{Mini-workshop on Lie groupoids, Lie algebroids and differentiable stacks} (June 11–14, 2024, Universit\`{a} degli Studi di Napoli Federico II), and a mini-course taught in April 2025 at the Instituto de Matem\'{a}tica e Estat\'{i}stica, Universidade de S\~{a}o Paulo (IME-USP). The author especially thanks Ivan Struchiner and the members of the IME-USP for their warm hospitality.

\section{Motivating examples and overview}\label{sec:MotivatingExamplesAndOverview}
Although certainly implicit in the work of Lie and Cartan, the first explicit connection between algebroids and classification problems was pointed out by Bryant in his work on Bochner-K\"{a}hler metrics \cite{Bryant2001}. Later, through many examples, he observed that many existence problems can be recast into a set of equations that resembled the structure of an algebroid but depends on extra ``free derivatives" \cite{Bryant2014}. In some special cases, the free derivatives disappear after a process called prolongation, and the set of equations defines a Lie algebroid \cite{Bryant2001, Bryant2014, FernandesStruchiner2019}. The integration theory for Lie algebroids leads to a beautical classification of the \emph{global solutions} of the geometric problem \cite{FernandesStruchiner2014, FernandesStruchiner2019, FernandesStruchiner2021}. Developing the Lie theory for generic geometric problems is the main motivation to introduce relative algebroids. In this section, we will discuss several examples where a geometric problem leads to a Lie algebroid or something that resembles an algebroid. Then we give an overview of how exploiting the Lie theory of this algebroid leads to insights into the geometric problem.

\subsection{Geometric structures and coframes}

A \textbf{coframe} on a $n$-dimensional manifold $P$ is a trivialization of the tangent bundle  $\theta\: TP\xrightarrow{\sim}\RR^n$. We denote by $(P, \theta)$ the manifold with coframe. The components of $\theta = (\theta^i)$ with respect to the standard basis on $\RR^n$ form a global frame of the cotangent bundle $T^*P$ (hence the name). Because of this, it is possible to express the exterior derivative on $P$ entirely in terms of the coframe:
\[
	\d \theta^i = -\frac{1}{2} c_{jk}^{i} \theta^j\wedge\theta^k.
\]
The functions $c^i_{jk}\in C^\infty(P)$ are called the \textbf{structure functions} of the coframe. Likewise, the components of the differential of a function $f\in C^\infty(P)$ in terms of $\theta$,
\[
	\d f= (\del_{\theta^i}f) \theta^i,
\] 
define the \textbf{coframe derivatives} of $f$. The vector fields $(\del_{\theta^i})$ constitute the frame of $TP$ dual to $(\theta^i)$. 

Coframes themselves already describe a wide variety of geometric structures, though not always directly on the manifold itself.

\begin{example}[Lie groups]\label{ex:IntroLieGroups}
	Let $G$ be a Lie groupoid with Lie algebra $\mf{g}$. The (left) Maurer-Cartan form $\theta_{\mathrm{MC}}$ is defined as
	\[
		\theta_{\mathrm{MC}}(X_g) = T_g L_{g^{-1}}(X_g)\in T_e G\cong \mf{g}
	\]
	for $X_g\in T_gG$, where $L_g\: G\to G$ is left-translation. After choosing a basis $(e_i)$ of $\mf{g}$, the components of the Maurer-Cartan form constitute a coframe on $G$. In this case, the structure functions
	\[
		\d \theta_{\mathrm{MC}}^i = -\frac{1}{2} c^i_{jk} \theta_{\mathrm{MC}}^j\wedge \theta_{\mathrm{MC}}^k
	\]
	are \emph{constant}! They correspond to structure constants of the Lie algebra $\mf{g}$ of $G$ with respect to the basis $(e_i)$:
	\[
		[e_j, e_k] = c^i_{jk}e_i.
	\]
	In fact, any manifold with coframe $(P, \theta)$ whose structure functions are constant is (locally) equivalent to the Maurer-Cartan form of some Lie group! We conclude that the local (resp. global) classification problem for coframes with constant structure functions is that of local (resp. global) Lie groups!
\end{example}

\begin{example}[Riemannian manifolds]\label{ex:IntroRiemannianManifolds}
	Let $(M, g)$ be an $n$-dimensional Riemannian manifold. While $M$ itself may not admit any coframe, the Riemannian structure $g$ is still captured by a coframe on a different space. The orthonormal frame bundle
	\[
		P:= \left\{ e\:(\RR^n, \langle\cdot, \cdot \rangle)\xrightarrow{\sim} (T_xM, g_x) \ | \  x\in M\right\}
	\]	
	consists of pointwise isomorphisms from $\RR^n$ with the standard inner product $\langle\cdot, \cdot\rangle$ to $(T_xM, g_x)$. It forms a principal $O(n)$-bundle $\pi\: P\to M$. The \textbf{tautological form} $\theta\in \Omega^1(P; \RR^n)$ on $P$ is defined to be
	\[
		\theta_e (X) = e^{-1}\left (T_e \pi (X)\right)\in \RR^n,
	\]
	for $X\in T_e P$. The tautological form itself is not a coframe because it has kernel equal to $\ker T\pi$, but for Riemannian manifolds, the Levi-Civita connection one-form $\omega\in \Omega^1(P; \mf{o}(n))$ is exactly complementary to $\theta$, so we obtain (after choosing a basis of $\mf{o}(n)$) a coframe
	\[
		\theta\oplus \omega \: TP \xrightarrow{\sim} \RR^n\oplus \mf{o}(n)
	\]
	canonically associated to the Riemannian manifold $(M, g)$.  The structure functions of $\theta\oplus \omega$ can be expressed basis-indepedently, and are known as Cartan's structure equations
	\begin{equation}\label{eq:CartansStructureEquationsRiemannianManifolds}
		\left\{
			\begin{aligned}
				&\d \theta = - \omega \wedge \theta, \\
				&\d \omega = R(\theta \wedge\theta) - \omega \wedge \omega,
			\end{aligned}
		\right.
	\end{equation}
	where $R\: P\to \Hom(\wedge^2\RR^n, \mf{o}(n))$ becomes the classical Riemann curvature tensor after passing to the associated bundles on $M$.
	
	Conversely, if $P$ is any manifold with a coframe $\theta\oplus \omega\: TP\to \RR^n \oplus \mf{o}(n)$ satisfying Cartan's Structure Equations (\ref{eq:CartansStructureEquationsRiemannianManifolds}), then it is \emph{locally} isomorphic to the orthonormal frame bundle of some Riemannian manifold.
\end{example}

\begin{example}[$G$-structures with connection]
	The reader familiar with $G$-structures might have recognized that the previous example generalizes to $G$-structures with a specified connection. The bundle of $G$-frames $\pi\:P\to M$ of a $G$-structure on $M$, together with a $G$-connection one-form $\omega$, carries a canonical coframe 
    \[
        \theta\oplus \omega\: TP\xrightarrow{\sim} \RR^n\oplus \mf{g}
    \]
    satisfying Cartan's structure equations
    \[
        \left\{
        \begin{aligned}
            &\d \theta = T(\theta \wedge\theta) - \omega \wedge \theta,\\
            &\d \omega = R(\theta\wedge\theta) - \omega\wedge\omega,
        \end{aligned}
        \right.
    \]
    where $T$ is the torsion and $R$ is the curvature of the connection, when passing to the associated bundles over $M$.
\end{example}

\subsection{Appearances of algebroids in classification problems for geometric structures}\label{sec:AlgebroidsInClassificationProblems}

Restricting to classification problems for Riemannian surfaces keeps the technical computations accessible while remaining conceptually clarifying. 

Let $(M, g)$ be a metric surface. After identifying $\mf{o}(2)\cong \RR\left(\begin{smallmatrix} 0 & -1 \\ 1 & 0 \end{smallmatrix}\right)$, we obtain a coframe on the orthonormal frame bundle $P$, with components
\[
    \theta = (\theta^1, \theta^2), \quad \omega = \theta^3.
\] 
where $\theta$ is the tautological form and $\omega$ the Levi-Civita connection (Example \ref{ex:IntroRiemannianManifolds}). Cartan's Structure Equations (\ref{eq:CartansStructureEquationsRiemannianManifolds}) expressed in these components read
\begin{equation}\label{eq:CartanStructureEquationsSurfaces}
    \left\{
        \begin{aligned}
            \d \theta^1 &= - \theta^3 \wedge \theta^2, \\
            \d \theta^2 &= \theta^3 \wedge \theta^1, \\ 
            \d \theta^3 &= K \theta^1\wedge \theta^2,
        \end{aligned}
    \right.
\end{equation}
where $K$ is the scalar curvature.

\begin{example}[Space forms]\label{ex:IntroSurfacesWithConstantCurvature}
	Suppose we are trying to classify Riemmanian surfaces with constant curvature. The condition that the scalar curvature is constant can be expressed through the equation $\d K = 0$, which can be adjoined to Equation (\ref{eq:CartanStructureEquationsSurfaces}) to obtain 
	\begin{equation}\label{eq:ClassificationConstantCurvature}
		\left\{
		\begin{aligned}
			\d \theta^1 &= - \theta^3 \wedge \theta^2, \\
			\d \theta^2 &= \theta^3 \wedge \theta^1, \\ 
			\d \theta^3 &= K \theta^1\wedge \theta^2, \\
			\d K 	&= 0.
		\end{aligned}
		\right.
	\end{equation}
	These are the equations that behind the classification problem for metric surfaces with constant curvature.

    A \emph{solution} is a manifold with coframe $(P, \theta^1, \theta^2, \theta^3)$ and a function $K\in C^\infty(P)$ such that the coframe satisfies Equation (\ref{eq:ClassificationConstantCurvature}).
    
	This notion of solution leads to believe that the essential structure behind Equation (\ref{eq:ClassificationConstantCurvature}) is obtained by \emph{forgetting} about the underlying manifold $P$ and interpreting the covectors and variables appearing in Equation (\ref{eq:ClassificationConstantCurvature}) as independent objects.\footnote{This is akin to treating the derivatives of a function to independent coordinates on a jet space.} In this case, the structure behind Equation (\ref{eq:ClassificationConstantCurvature}) consists of
	\begin{itemize}
		\item a parameter space $\RR = \{K\}$,
		\item a vector bundle $\underline{\RR^3}\to \RR$, with global frame $(e_1, e_2, e_3)$ dual to $(\theta^1, \theta^2, \theta^3)$,
		\item a degree-1 derivation $``\d" = \D \: \Gamma(\wedge^\bullet (\underline{\RR^3})^*)\to \Gamma(\wedge^{\bullet+1} (\underline{\RR^3})^*)$ satisfying $\D^2 =0$. 
	\end{itemize}
    Here, the underline in $\underline{\RR^3}\to \RR$ indicates the trivial vector bundle with fiber $\RR^3$. 

    A degree-1 derivation on a vector bundle that squares to zero is exactly what defines a Lie algebroid! 
\end{example}

\begin{example}[Metrics of Hessian type]
Following \cite[Section 5.2]{Bryant2014}, a metric surface $(M, g)$ is of Hessian type when the Hessian of the scalar curvature satisfies the equation
\begin{equation}\label{eq:HessianSurface}
    \operatorname{Hess}_g(K) = a(K)g + b(K) \d K^2
\end{equation}
for some functions $a, b$ of one variable. To rephrase this equation in terms of the orthonormal frame bundle, we write
\begin{align*}
    \d K &= K_1 \theta^1 + K_2\theta^2, & \d K_1 &= K_{1i}\theta^i, & \d K_{2} &= K_{2i} \theta^i.
\end{align*}
From $\d^2 K=0$ it follows that $K_{12}= K_{21}$, $K_{13} = - K_2$ and $K_{23}= K_1$. The Hessian is expressed as $\operatorname{Hess}_g(K) = K_{11} \left(\theta^1\right)^2 + 2K_{12} \theta^1\theta^2 + K_{22}\left(\theta^2\right)^2$. Equation (\ref{eq:HessianSurface}) implies that
\begin{equation}\label{eq:ClassificationHessianSurfaces}
    \left\{
        \begin{aligned}
            & \d K = K_1 \theta^1 + K_2\theta^2,\\
            & \d K_1 = (a(K) + b(K)K_1^2)\theta^1 + b(K)K_1K_2 \theta^2 - K_2 \theta^3,\\
            & \d K_2 = b(K)K_1K_2\theta^1 + (a(K) + b(K)K_2^2) \theta^2 + K_1 \theta^3.
        \end{aligned}
    \right. 
\end{equation}
Adding these equations to Equation (\ref{eq:CartanStructureEquationsSurfaces}) yields the full set of equations governing this problem.

Again, the next step is to forget about the underlying manifold. Interpreting Equations (\ref{eq:CartanStructureEquationsSurfaces}) and (\ref{eq:ClassificationHessianSurfaces}) 
\begin{itemize}
    \item a parameter space $\RR^3 = \{ K_1, K_2, K\}$,
    \item a vector bundle $\underline{\RR^3}\to \RR$ with global frame $(e_1, e_2, e_3)$ dual to $(\theta^1, \theta^2, \theta^3)$,
    \item a degree-1 derivation $``\d" = \D \: \Gamma(\wedge^\bullet (\underline{\RR^3})^*)\to \Gamma(\wedge^{\bullet+1} (\underline{\RR^3})^*)$.
\end{itemize}
    However, unlike in Example \ref{ex:IntroSurfacesWithConstantCurvature}, the condition $\D^2=0$ is not automatically satisfied. This time, $\D^2=0$ implies that $(a'(K) - a(K) b(K) + K)K_i =0$. Either we have to restrict to the locus $K_1=K_2=0$, which leads to the same algebroid as in Example \ref{ex:IntroSurfacesWithConstantCurvature}, or we assume that $a'(K) = a(K) b(K)-K$, in which case we obtain a Lie algebroid over the base $\RR^3$. 
\end{example}

\begin{example}[Extremal K\"ahler surfaces] Our next example is taken from \cite{FernandesStruchiner2019}, which is based on \cite{Bryant2001}. This time, we are looking at \emph{K\"{a}hler} surfaces. Since $\mf{u}(1)\cong \mf{o}(2)$, we can still use Equation (\ref{eq:CartanStructureEquationsSurfaces}) and change the classification problem to K\"{a}hler surfaces.

A K\"{a}hler surface $(M, g, \omega)$ is \textbf{extremal} when the Hamiltonian vector field of the scalar curvature is Killing. To derive the equations on the level of the (oriented) orthormal frame bundle, let $X_K$ be the Hamiltonian vector field associated to the scalar curvature $K$. Then the assumption that $X_K$ is Killing implies that its complete lift preserves the orthonormal frame bundle $\pi\:P\to M$ and the Levi-Civita connection. So, as complete lifts always preserve the tautological form, we obtain a vector field $\tilde{X}_K$ on $P$ that preserves the coframe $(\theta^1, \theta^2, \theta^3)$. This lift satisfies
\[
    \iota_{\tilde{X}_K} (\theta^1\wedge\theta^2) = \iota_{\tilde{X}}\pi^*\omega = \pi^*\d K = K_1\theta^1 + K_2\theta^2, 
\]
which determines two components of $\tilde{X}_K$. The third component $\theta^3(\tilde{X}_K)$ we denote by $U$. 

Then, since $\LL_{\tilde{X}_K}\theta =0$, we find
\begin{equation}
    \left\{
    \begin{aligned}
        & \d K = K_1 \theta^1 + K_2 \theta^2,\\
        & \d K_1 = - \d \iota_{\tilde{X}_K} \theta^2 = \iota_{\tilde{X}_K}\d \theta^2 = U\theta^1 - K_2 \theta^3,\\
        & \d K_2 = \d \iota_{\tilde{X}_K}\theta^1 = -\iota_{\tilde{X}_K} \d \theta^1 = U \theta^2 + K_1\theta^3,\\
        & \d U = \d \iota_{\tilde{X}_K} \theta^3 = -\iota_{\tilde{X}_K}\d\theta^3= - K \d K. 
    \end{aligned}
    \right.
\end{equation}
These are the equations governing the classification of extremal K\"{a}hler surfaces. This time, the equatios define a Lie algebroid $\underline{\RR^3}\to \RR^4= \{U, K_1, K_2, K\}$.
\end{example}

The appearance of an actual Lie algebroid is rather special. For this to happen, the system of equations determining the structure has to ``close up" in order to define a degree-1 derivation on a vector bundle over a finite-dimensional space. Generically, an equation on the curvature does not lead to such a closed system, and the moduli space of local solutions is infinite-dimensional. Therefore, studying such problems from the perspective of Lie theory asks for a broader notion of algebroid. Let us illustrate this behaviour with an example, many more can be found in \cite{Bryant2014}.

\begin{example}[Surfaces with $|\nabla K| = 1$]\label{ex:IntroSurfacesWith|nablaK|=1}
    Suppose that we want to study metric surfaces with scalar curvature satisfying the diffeomorphism-invariant equation $|\nabla K| = 1$ (as per \cite[Section 5.1]{Bryant2014}). On the level of the orthonormal frame bundle, this equation translates to
    \[
        \d K = K_1 \theta^1 + K_2 \theta^2, \quad K_1^2 + K_2^2 = 1.
    \]
    We can parametrize the derivative $\d K$ by the coordinate $\varphi$ on the circle $S^1$ to obtain the set of equations governing this classification problem:
    \begin{equation}\label{eq:Classification|nablaK|=1}
		\left\{
		\begin{aligned}
			\d \theta^1 &= - \theta^3 \wedge \theta^2, \\
			\d \theta^2 &= \theta^3 \wedge \theta^1, \\ 
			\d \theta^3 &= K \theta^1\wedge \theta^2, \\
			\d K 	&= (\cos \varphi)\theta^1 + (\sin \varphi)\theta^2.
		\end{aligned}
		\right.
	\end{equation}
    The function $\varphi$ is left free.

    A solution to Equation (\ref{eq:Classification|nablaK|=1}) is a manifold with coframe $(P, \theta^1, \theta^2, \theta^3)$ and a map $(K, \varphi)\: P\to \RR\times S^1$ that together satisfy Equations (\ref{eq:Classification|nablaK|=1}).

    This time, the structure behind the classifying equations is much more elusive. Treating the variables as independent objects, we see that we have
    \begin{itemize}
        \item a parameter space $\RR = \{K\}$,
        \item a bundle of free derivatives over the parameter space $p\: S^1\times \RR\to \RR$, where $S^1\times \RR$ has coordinate $(\varphi, K)$, and $p$ is the projection,
        \item a vector bundle $A = \underline{\RR^3}\to \RR$ spanned by a frame $(e_1, e_2, e_3)$ dual to $(\theta^1, \theta^2, \theta^3)$,
        \item a derivation $``d" = \D$ that is defined on $\Gamma(\wedge^\bullet A^*)$. However, when we apply $\D$ to a form $\alpha \in \Gamma(\wedge^\bullet A^*)$, the result $\D \alpha$ depends also on the free derivative $\varphi$. Thus, if $B = p^*A\to S^1\times \RR$ is the pullback bundle, then $\D$ becomes a map
        \[
            \D\: \Gamma(\wedge^\bullet A^*) \to \Gamma(\wedge^\bullet B^*),
        \]
        which is a derivation \emph{relative to $p$}.
    \end{itemize}

    One could argue that Equation (\ref{eq:Classification|nablaK|=1}) is not a complete set of equations. The derivatives of $\varphi$ are left unspecified, even though they are not entirely unrestricted. The set of equations is still subject ``$\d^2=0$", which indeed gives restrictions on $\d \varphi$:
    \begin{align*}
        0 &= \d^2 K = -(\sin \varphi)\d \varphi \wedge \theta^1 + (\cos \varphi) \d \theta^1 + (\cos \varphi)\d \varphi \wedge \theta^2 + (\sin \varphi) \d \theta^2\\
        &= (\d \varphi - \theta^3) \wedge ( - (\sin \varphi) \theta^1 + (\cos \varphi) \theta^2).
    \end{align*}
    It follows that $\d \varphi = \theta^3 + c_1 ( - (\sin \varphi) \theta^1 + (\cos \varphi) \theta^2)$ for a \emph{new} free variable $c_1$. Adding it to the set of equations would still make the resulting derivation relative. The procedure of obtaining new free derivatives by restraining the old free derivatives is called \emph{prolongation}.
\end{example}

\subsection{From algebroids to classification}

As illustrated by the examples in the previous section, there are generically two types of classification problems: those of \emph{finite type}, where the system closes up after finitely many prolongations, and those of \emph{infinite type}, where at each prolongation new free derivatives are introduced.

\subsubsection{Finite type classification problems}
    The theory developed by Fernandes and Struchiner in \cite{FernandesStruchiner2014, FernandesStruchiner2021, FernandesStruchiner2019} exploits beautifully the Lie theory behind finite-type classification problems.

    Suppose that we are given a Lie algebroid $(A, \D_A)\to M$ that comes from a classification problem (as in the examples from the previous section). This means that the \emph{solutions} to the geometric problem are manifolds $(P, \theta, r)$ with a Lie algebroid morphism
    \[
        \xymatrix{
            (TP, \d) \ar[r]^{\theta} \ar[d] & (A, \D_A) \ar[d] \\
            P \ar[r]^{r} & M
        }
    \]
    that is fiberwise an isomorphism.

    The theory of integration of Lie algebroids leads to a complete understanding of the global solutions. Suppose that $\GG\toto M$ is a Lie groupoid integrating $(A, \D_A)$. If we fix $x\in M$, then the source fiber with the Maurer-Cartan form $(s^{-1}(x), \theta_{\mathrm{MC}}, t)$ is a \emph{complete} solution, in the sense that it can not be extended to a strictly larger solution. 

    Moreover, every complete \textit{simply connected} solution $(P, \theta, r)$ is isomorphic to a source fiber of the \emph{source-simply connected} integration of $(A, \D_A)$!

\begin{remark}
    The story for finite type problems as presented here is oversimplified for the sake of exposition. There is a lot of subtlety involved especially related to the integration problem and the role of the structure group of the $G$-structure. This subtlety is fully understood through the framework of $G$-structure algebroids and $G$-structure groupoids \cite{FernandesStruchiner2019}.
\end{remark}

\subsubsection{Infinite type classification problems}
    Generically, classification problems are of infinite type. As in example \ref{eq:Classification|nablaK|=1}, the equations for such problems define \textbf{derivation} $\D$ on a vector bundle $A\to N$ \textbf{relative to a submersion} $p\: M\to N$. Denoting $B:= p^*A$, we like to depict this structure by
    \[
        \xymatrix{
            B \ar[r] \ar[d] & A \ar[d] \ar@/^1pc/@{-->}[l]^{\D} \\
            M \ar[r]^p & N
        }
    \]
    where the dashed arrow is not a map but indicates that the derivation is defined on $\Omega^\bullet(A):= \Gamma(\wedge^\bullet A^*)$ but takes values in $\Omega^{\bullet+1}(B)$. 

    In this case, a \textbf{solution} (or \textbf{realization}) is a manifold with a morphism of relative algebroids $(P, \theta, r)$ that is fiberwise an isomorphism:
    \[
        \xymatrix{
            (TP, \d) \ar[d] \ar[r]^\theta & (B, \D) \ar[d] \\
            P \ar[r]^r & M, 
        }
    \]
    That is, $\theta^*\circ \D = \d \circ \theta^*\circ p^*$ on $\Omega^\bullet(A)$. Such a morphism ``realizes" the derivation $\D$ as the de Rham differential on the manifold $P$.

    Successive prolongations (if they exist) gives rise to a sequence of relative algebroids
    \[
    \xymatrix{
    {\left(B^{(\infty)}, \D^{(\infty)}\right)\:\ldots} \ar[r] \ar@<-12pt>[d] & B^{(k)} \ar[d] \ar[r] &  \ldots \ar[r] \ar@{-->}@/^1pc/[l]^{\D^{(k)}} & B^{(1)} \ar[d] \ar[r] & B \ar@{-->}@/^1pc/[l]^{\D^{(1)}} \ar[r] \ar[d] & A \ar[d] \ar@{-->}@/^1pc/[l]^{\D} \\
    M^{(\infty)}\:\ldots \ar[r]                                            & M^{(k)} \ar[r]_{p_k}    & \ldots \ar[r]                                                               & M^{(1)} \ar[r]_{p_1} & M \ar[r]_{p}    & N
    }               
    \]
    where
    \begin{itemize}
        \item each vector bundle is the pullback of the previous one: $B^{(k+1)} = p_{k+1}^*B^{(k)}$,
        \item each derivation \emph{extends} the previous one: $p_{k+1}^* \circ \D^{(k)} = \D^{(k+1)}\circ p_k^*$,
        \item consecutive derivations square to zero: $\D^{(k+1)} \circ \D^k = 0$
    \end{itemize}
    Together, they assemble into an honest derivation $\D^{(\infty)}$ that squares to zero on the profinite sections of $B^{(\infty)}\to M^{(\infty)}$. Therefore, the resulting object $(B^{(\infty)}, \D^{(\infty)})$ is a \emph{profinite Lie algebroid}.

    Guided by the finite type theory, the analysis of the profinite algebroid related to a classification problem already leads to a lot of insights.
    \begin{itemize}
        \item The space $M^{(\infty)}$ can be interpreted as the space of (pointed) \emph{formal} solutions modulo formal diffeomorphism.
        \item Since the vector bundle $B^{(\infty)}$ has finite rank, the isotropy Lie algebras are finite dimensional Lie algebras, and correspond to infinitesimal symmetries of the solutions.
        \item The profinite Lie algebroid can detect when two realizations passing though the same point are equivalent.
    \end{itemize}
    
    The profinite Lie algebroid behind a classification problem also leads to a lot of natural but hard questions. First of all, what does an integration of $(B^{(\infty)}, \D^{(\infty)})$ mean? We believe that it can not be a ``profinite Lie groupoid", but that it must be an object of entirely different nature. 

    If the profinite Lie algebroid $(B^{(\infty)}, \D^{(\infty)})$ is \emph{real analytic}, and if we restrict to \emph{real analytic solutions}, then much of the finite dimensional geometry carries over. For instance, the locus $M^\omega \subset M^{(\infty)}$ of points in the image of a real analytic solutions is partitioned into finite-dimensional leaves, and on each leaf, the algebroid restricts to a finite-dimensional Lie algebroid!

    This divides the problem of finding the complete real analytic solutions into two steps: analytic continuation of the analytic leaves of the algebroid, and Lie algebroid integration in the isotropy direction. The first step is generally hard, but the second is very well understood!

We hope this overview inspires the reader to dive deeper into the theory of relative algebroids. The rest of these notes provide an introduction to relative algebroids and their prolongation theory.

\section{Relative derivations}\label{sec:RelativeDerivations}

In this section, we investigate the structure of derivations relative to a map or a smooth submersion. It is possible to consider derivations relative to arbitrary bundle maps or even foliations, but for our purposes we stick to submersions, and refer to \cite{FernandesSmilde2025} to these more general notions.

\subsection{Derivations, brackets and algebroids}\label{sec:DerivationsBracketsAndAlgebroids}
We start by recalling the structure of ordinary derivations on vector bundles, as developed in for instance \cite{CrainicMoerdijk08}. 

Let $E\to M$ be a vector bundle. We denote by $\Omega^\bullet(E) := \Gamma(\wedge^\bullet E^*)$ the exterior algebra of $E$-forms. A \textbf{$k$-derivation} on $E$ is a map
\[
    \D\: \Omega^\bullet(E)\to \Omega^{\bullet + k}(E)
\]
satisfying the Leibniz rule:
\[
    \D(\alpha\wedge \beta) = (\D\alpha)\wedge \beta + (-1)^{|\alpha| k} \alpha \wedge (\D \beta). 
\]
for all homogenous $\alpha, \beta \in \Omega^\bullet(E)$. We denote by $\Der^k(E)$ the set of $k$-derivations on $E$. 

Dual to a $k$-derivation, there is the notion of a $k$-bracket on $E$, which consist of a bracket and an anchor
\[
    [\cdot, \ldots, \cdot]\: \wedge^{k+1}\Gamma(E)\to \Gamma(E), \quad \rho\: \wedge^k E\to TM
\]
subject to the Leibniz rule:
\[
    [v_1, \dots, fv_{k+1}] = f[v_1, \dots, v_{k+1}] + \LL_{\rho(v_1, \dots, v_k)}(f) v_{k+1}
\]
for $v_1, \dots, v_{k+1}\in \Gamma(E)$ and $f\in C^\infty(M)$. 

Every $k$-derivation determines and is determined by a $k$-bracket through the Koszul formula, presented here only on $\Omega^0(E)$ and $\Omega^1(E)$:
\begin{equation}\label{eq:KoszulFormula}
    \begin{aligned}
    (\D f)(v_1, \dots, v_k) &= (\d f)(\rho(v_1, \dots, v_k)),\\
    (\D \alpha) (v_0, \ldots, v_k) &= \sum_{i=0}^k \LL_{\rho(v_0, \dots, \widehat{v_i}, \dots, v_k)} (\alpha(v_i)) - \alpha([v_0, \dots, v_k]).
    \end{aligned}
\end{equation}
for $f\in C^\infty(E)$, $\alpha\in \Omega^1(E)$ and $v_0, \dots, v_k\in \Gamma(E)$. 

\begin{example}[Exterior derivative]
    On $E=TM$, the tangent bundle of a manifold $M$, the de Rham differential $\d$ is a degree-1 derivation. The 1-bracket dual to $\d$ is the Lie bracket of vector fields $[\cdot, \cdot]$ and anchor $\rho = \id_{TM}$. 
\end{example}

\begin{example}[Lie algebras]
    If $\mf{g}$ is a Lie algebra, the Chevalley-Eilenberg differential $\d_{\mathrm{CE}}$ is the derivation dual to the Lie bracket $[\cdot, \cdot]_{\mf{g}}$ on $\mf{g}$. Explicitly, 
    \[
        \d_{\mathrm{CE}}\alpha(v, w) = - \alpha([v, w]_{\mf{g}})
    \]
    for $\alpha\in \mf{g}^*$ and $v, w\in \mf{g}$.
\end{example}

\begin{example}[Linear vector fields]
    A 0-derivation $\D\: \Omega^\bullet(E)\to \Omega^\bullet(E)$ corresponds to a linear vector field $X\in \mf{X}(E)$ whose flow $\Phi_t$ satisfies
    \[
    \left\{
        \begin{aligned}
            \frac{\d}{\d t} \Phi^*_t(\alpha) 
            &= \Phi^*_t(\D \alpha), \\
            \Phi_0 &= \id_E,
        \end{aligned}
    \right.
    \]
    for $\alpha\in\Omega^1(E)$. 
\end{example}

\begin{example}[Square of a 1-derivation]
    If $\D$ is a 1-derivation on $E$, then $\D^2$ is a 2-derivation on $E$. If $[\cdot, \cdot]$ is the bracket dual to $\D$, then the bracket dual to $\D^2$ is the Jacobiator of $[\cdot, \cdot]$. 
\end{example}

\begin{example}
    It is instructive to write a derivation $\D$ in local coordinates. For simplicity, we assume that $\D$ is a 1-derivation on $E$---similar formulas can be derived for arbitrary $k$-derivations. Choose coordinates $(x^\mu)$ on $M$ and a local frame $(\theta^i)$ of $E^*$. Then $\D$ is determined by
    \begin{equation}\label{eq:DerivationInCoordinates}
        \left\{
            \begin{aligned}
                &\D \theta^i = -\frac{1}{2}c^i_{jk} \theta^j\wedge\theta^k, \\
                &\D x^\mu = \rho^\mu_i \theta^i.
            \end{aligned}
        \right.
    \end{equation}  
    The functions $c^{i}_{jk}, \rho^\mu_i\in C^\infty(M)$ are the \textbf{structure functions} of $\D$ with respect to these coordinates.  
    
    The bracket dual to $\D$ in terms of the structure functions is given by
    \begin{align*}
            [e_j, e_k] &= c_{jk}^i e_i, & \rho(e_i) = \rho^\mu_i \frac{\del}{\del x^\mu},
    \end{align*}
    where $(e_i)$ is the frame of $E$ to $(\theta^i)$. 
\end{example}

\subsubsection{The structure of derivations}

Let $E\to M$ be a vector bundle. The space of $k$-derivations on $E$ is a $C^\infty(M)$ module, meaning that if $\D\in \Der^k(E)$, then also $f\D \in \Der^k(E)$ for any $f\in C^\infty(M)$. Equation (\ref{eq:DerivationInCoordinates}) shows that the $C^\infty(M)$-module $\Der^k(E)$ is locally finitely generated, so by the Serre-Swan theorem, there exists a vector bundle 
\[
    \DD^k_E\to M
\]
with $\Gamma(\DD^k_E) \cong \Der^k(E)$. The bundle $\DD^k_E$ is called the \textbf{bundle of (pointwise) $k$-derivations} on $E$. We will see later in Example \ref{ex:PointwiseDerivations} that we can interpret elements $\D_x\in \DD^1_E$ as derivations relative to the inclusion $\{x\}\hookrightarrow M$.

Finally, the \textbf{symbol} $\sigma(\D)\: \wedge^kE\to M$ of a derivation $\D\in \Der^k(E)$ is precisely the anchor of the corresponding $k$-bracket, defined through the Koszul formula (\ref{eq:KoszulFormula}). The symbol defines a bundle map $\sigma\: \DD^k_E\to \Hom(\wedge^kE, TM)$ that fits into the \textbf{symbol short exact sequence}
\begin{equation}\label{eq:SymbolSES}
    \xymatrix{
    0 \ar[r] & \Hom(\wedge^{k+1}E, E) \ar[r] & \DD^k_E \ar^-\sigma[r] & \Hom(\wedge^k E, TM) \ar[r] & 0.
    }
\end{equation}

\subsubsection{Lie algebroids}
A \textbf{Lie algebroid} $(A, \D)$ is a vector bundle $A\to M$ with a 1-derivation $\D \in \Der^1(A)$ such that $\D^2 = 0$. Equivalently, a Lie algebroid $(A, [\cdot, \cdot], \rho)$ consists of a vector bundle $A\to M$ together with a 1-bracket $[\cdot, \cdot]\: \wedge^2\Gamma(A)\to \Gamma(A)$ with anchor $\rho\: A\to TM$ such that $[\cdot, \cdot]$ satisfies the \emph{Jacobi identity}.

\begin{example}[Lie algebroid of surfaces with constant curvature]
In the classification of surfaces with constant curvature (Example \ref{ex:IntroSurfacesWithConstantCurvature}), the equations behind the classification problem (Equation (\ref{eq:ClassificationConstantCurvature})) define a Lie algebroid on the trivial bundle $(\underline{\RR^3}, \D)\to \RR$, where the base space $\RR$ has coordinate $\{K\}$. From the equation $\D K =0$, it follows from the Koszul formula (\ref{eq:KoszulFormula}) that the anchor of this Lie algebroid vanishes. So $(\underline{\RR^3}, \D)$ is a bundle of Lie algebras, with isotropy given by
\[
    \begin{cases}
        \mf{o}(3), &\mbox{when $K>0$},\\
        \mf{o}(2)\ltimes \RR^2, &\mbox{when $K=0$}, \\
        \mf{o}(2, 1), & \mbox{when $K<0$}.
    \end{cases}
\]
\end{example}

\subsection{Relative derivations and brackets}

Let $A\to N$ be a vector bundle and $p\: M\to N$ be any smooth map (typically a submersion). We denote by $B:= p^*A$ the pullback bundle. 
\begin{definition}
    A \textbf{$k$-derivation} on $A$ \textbf{relative to $p$} is a map
    \[  
        \D \: \Omega^\bullet(A) \to \Omega^{\bullet + 1}(B)
    \]
    subject to the Leibniz rule:
    \[
        \D(\alpha\wedge\beta) = (\D \alpha) \wedge (p^*\beta) + (-1)^{|\alpha| k} (p^*\alpha) \wedge (\D \beta) 
    \]
    for homogenous $\alpha, \beta\in \Omega^\bullet(A)$. We denote by $\Der^k(A, p)$ the space of $k$-derivations on $A$ relative to $p$.
\end{definition}    

An ordinary derivation is a derivation relative to the identity. Therefore, the examples in Section \ref{sec:DerivationsBracketsAndAlgebroids} are in particular relative derivations. 

Like for ordinary derivations, there is a notion of a relative bracket dual to that of a relative derivation.

\begin{definition}\label{def:RelativeBrackets}
    A $k$-bracket relative to $p$ consists of a bracket and an anchor
    \begin{align*}
        [\cdot, \dots, \cdot]&\:\wedge^{k+1}\Gamma(A)\to \Gamma(B), & \rho&\: \wedge^k B\to p^*TN,
    \end{align*}
    such that
    \[
        [v_1, \dots, fv_{k+1}] = (p^*f)[v_1, \dots, v_{k+1}] + \LL_{\rho(p^*v_1, \dots, p^*v_k)}(f) v_{k+1}
    \]
    for $v_1, \dots, v_{k+1}\in \Gamma(A)$ and $f\in C^\infty(N)$.
\end{definition}

Note that the derivate $\LL_{\rho(p^*v_1, \dots, p^*v_k)}(f) = \langle \d f, \rho(p^*v_1, \dots, p^*v_k)\rangle$ makes sense \emph{pointwise} on $M$ because $\rho$ takes values in $p^*TN$, so the result after evaluation is a function on $M$. 

As before, relative derivations and relative brackets are in duality through the relative version of the Koszul formula:
\begin{equation}\label{eq:RelativeKoszulFormula}
    \begin{aligned}
    (\D f)(p^*v_1, \dots, p^*v_k) &= (\d f)(\rho(p^*v_1, \dots, p^*v_k)),\\
    (\D \alpha) (p^*v_0, \ldots, p^*v_k) &= \sum_{i=0}^k \LL_{\rho(p^*v_0, \dots, \widehat{p^*v_i}, \dots, p^*v_k)} (\alpha(v_i)) - \alpha([v_0, \dots, v_k]).
    \end{aligned}
\end{equation}
for $f\in C^\infty(M)$, $\alpha\in \Omega^\bullet(A)$ and $v_0, \dots, v_k\in \Gamma(A)$. 

Relative derivations appear in a wide variety of contexts, but we leave the concrete examples to Section \ref{sec:RelativeAlgebroids}. The following two examples are more abstract but fundamental later theory and computations.

\begin{example}[Pointwise derivations]\label{ex:PointwiseDerivations}
    Let $A\to N$ be a vector bundle and $x\in N$. A $k$-derivation relative to the inclusion $i_x\:\{x\}\to N$ is a map
    \[
        \D_x\: \Omega^\bullet(A)\to \wedge^{\bullet+k}A^*_x,
    \]
    subject to the Leibniz rule. Every derivation $\D\in \Gamma(\DD^k_A)$ can be composed with evaluation at $x$ to obtain a derivation $\D_x$ relative to $i_x$. This derivation only depends on the value of $\D$, regarded as a section of $\DD^k_A$, at $x\in N$. Therefore, we can interpret:
    \begin{equation}\label{eq:PointwiseDerivations}
        \Der^k(A,i_x) \cong (\DD^k_A)_x\cong i_x^* \DD^k_A.
    \end{equation}
    This observation will be used frequently throughout these notes.
\end{example}

\begin{example}[Relative 1-derivations in coordinates]
Let $\D$ be a 1-derivation on $A$ relative to a submersion $p\: M\to N$. Choose submersion coordinates $(y^\varrho, x^\mu)$ for $p$ of $M$ and $N$, where $p$ becomes the projection $(y^\varrho, x^\mu)\mapsto x^\mu$, and a local frame $(\theta^i)$ of $A^*$. The derivation $\D$ is determined by
\begin{equation}\label{eq:StructureEquationsRelativeDerivations}
\left\{
    \begin{aligned}
        &\D \theta^i = -\frac{1}{2} c^i_{jk}  \theta^j\wedge\theta^k, \\
        &\D x^\mu = \rho^\mu_i \frac{\del}{\del x^\mu}.
    \end{aligned}
\right.
\end{equation}
As before, the coefficients $c^{i}_{jk}, \rho^\mu_i\in C^\infty(M)$ are called the \textbf{structure functions} of the relative derivation (with respect to the local frame and coordinates). The derivation acting on the fiber coordinates $(y^\varrho)$ is left unspecified. Together with the dependency of the structure functions on $(y^\varrho)$, this is exactly what makes it a relative derivation. 

The relative bracket $[\cdot, \cdot]$ and anchor $\rho$ dual to $\D$ is determined by the structure functions:
\begin{align*}
    [e_j, e_k] &= c^i_{jk}(y^\varrho, x^\mu) e_i,  & \rho(e_i) = \rho^\mu_i(y^\varrho, x^\mu) \frac{\del}{\del x^\mu},
\end{align*}
where $(e_i)$ is the frame of $A$ dual to $(\theta^i)$. 
\end{example}

\subsection{Structure of relative derivations}\label{subsec:StructureOfRelativeDerivations}

Let $A\to N$ be a vector bundle and $p\: M\to N$ a smooth map. The structure of the bundle of relative derivations is more intricate than that of ordinary derivations. We will see that besides a symbol short exact sequence, there is also a component in the direction of the vertical coordinates. 

The space of $k$-derivation $\Der^k(A, p)$ is this time a $C^\infty(M)$-module. Equation (\ref{eq:StructureEquationsRelativeDerivations}) shows that this module is locally finitely generated, and therefore, again by the Serre-Swan theorem, there is a vector bundle $\DD^k_{(A, p)}\to M$ such that $\Gamma(\DD^k_{(A, p)})\cong \Der^k(A, p)$. 

Every derivation $\D\in \Der^k(A, p)$ can be restricted to a point $x\in M$ to obtain a map
\[
    \D_x\: \Omega^\bullet(A)\to \wedge^{\bullet + k}A^*_{p(x)},
\]
which is a pointwise derivation on $A$ (Example \ref{ex:PointwiseDerivations}). This shows that at least fiberwise $(\DD^k_{(A, p)})_x\cong (\DD^k_A)_{p(x)}$. This holds true at the level of vector bundles: according to [\cite[Lemma 1.9]{FernandesSmilde2025}], there is a canonical identification $\DD^k_{(A, p)}\cong p^*\DD^k_A$. Consequently, there is an isomorphism
\begin{equation}
    \Der^k(A, p)\cong \Gamma(p^*\DD^k_A).
\end{equation}

We are still writing $B=p^*A$. Every derivation $\D\in \Der^k(B)$ can be precomposed with $p^*\: \Omega^\bullet(A)\to \Omega^\bullet(B)$ resulting in a  derivation $p_*(\D) := \D \circ p^*$ on $A$ relative to $p$. This assignment is $C^\infty(M)$-linear and therefore descends to a bundle map $p_*\: \DD^k_B\to p^*\DD^k_A$ which, when $p$ is a \emph{submersion}, fits into a short exact sequence
\[
    \xymatrix{
        0 \ar[r] & \Hom(\wedge^k B, \ker Tp) \ar[r] & \DD^k_B\ar[r]^{p_*} & p^*\DD^k_{A} \ar[r] & 0.
    }
\]

Every relative derivation $\D\in \Der^k(A, p)$, just like ordinary derivations, is also assigned a \textbf{symbol} $\sigma(\D)\in \Hom(\wedge^k B, p^*TN)$, which coincides with the anchor of the corresponding relative bracket, defined through the relative Koszul formula (\ref{eq:RelativeKoszulFormula}). The symbol is a natural assignment that induces a map of short exact sequences
\begin{equation*}
    \xymatrix@C=0.04\textwidth{
        {} & 0 \ar[r] \ar[d] & \Hom(\wedge^kB, \ker Tp) \ar[r]^-\sim \ar[d] & \Hom(\wedge^k B, \ker Tp) \ar[d] \ar[r]  &0 \\
        0 \ar[r] & \Hom(\wedge^{k+1}B, B) \ar[d]^{p_*} \ar[r] & \DD^k_B \ar[r]^-\sigma \ar[d]^{p_*} & \Hom(\wedge^kB, TM) \ar[d]^{p_*} \ar[r] & 0\\
        0 \ar[r] & p^*\Hom(\wedge^{k+1}A, A) \ar[r] & p^*\DD^k_A \ar[r]^-\sigma & p^*\Hom(\wedge^k A, TN)\ar[r] & 0.
    }
\end{equation*}
Since $p$ is a submersion, then the columns are also exact. This diagram describes the structure of relative derivations.

\section{Relative algebroids}\label{sec:RelativeAlgebroids}
Now that we have discussed the foundations of relative derivations, we are ready to introduce the concept of a relative algebroid and study examples. 

\begin{definition}
    An \textbf{algebroid relative to a submersion}, or \textbf{relative algebroid}, consists of 
    \begin{itemize}
        \item a vector bundle $A\to N$,
        \item a submersion $p\: M\to N$, 
        \item a 1-derivation $\D\in \Gamma(p^*\DD^1_A)$ relative to $p$. 
    \end{itemize}
    We denote the structure of a relative algebroid by $(A, p, \D)$, and the pullback bundle by $B= p^*A$. Diagramatically, we like to depict a relative algebroids as
    \[
        \xymatrix{
            B \ar[r] \ar[d] & A \ar@/^1pc/@{-->}[l]^{\D} \ar[d] \\
            M \ar[r]^p & N.
        }
    \]
\end{definition}
In submersion coordinates $p\: (y^\varrho, x^\mu)\mapsto x^\mu$ and a local frame $(\theta^i)$ of $A^*$, the relative algebroid is determined by the structure functions of the relative derivation $\D$, which according to Equation (\ref{eq:StructureEquationsRelativeDerivations}) look like
\begin{equation}\label{eq:RElativeAlgebroidInCoordinates}
    \left\{
        \begin{aligned}
            &\D \theta^i = - \frac{1}{2} c^i_{jk} (y^\varrho, x^\mu) \theta^j\wedge \theta^k,\\
            & \D x^\mu = \rho^\mu_i(y^\varrho, x^\mu) \theta^i.
        \end{aligned}
    \right.
\end{equation}

\begin{definition}
    A \textbf{realization} $(P, \theta, r)$ of a relative algebroid $(A, p, \D)$ is a manifold with a bundle map $(\theta, r)\: TP\to B$ covering $r\: P\to M$ such that 
    \begin{equation}\label{eq:RealizationCondition}
        \d \circ \theta^* \circ p^* = \theta^*\circ \D , \quad \mbox{on $\Omega^\bullet(A)$.}
    \end{equation}
\end{definition}
A realization ``realizes" the derivation $\D$ as the de Rham differential on the manifold $P$. In local coordinates, the realization form is given by components $(\theta^i)$ and the realization map $r = (b^\varrho, a^\mu)$. Equation (\ref{eq:RealizationCondition}) is then equivalent to
\begin{equation}\label{eq:BryantsEquations}
    \left\{ 
    \begin{aligned}
        &\d \theta^i = -\frac{1}{2}c^i_{jk}(b^\varrho, a^\mu)\theta^j\wedge\theta^k,\\
        &\d a^\mu = \rho^\mu_i(b^\varrho, a^\mu) \theta^i.
    \end{aligned}        
    \right.
\end{equation}
These equations appeared in \cite{Bryant2014} and therefore we call them \textbf{Bryant's equations}. It is clear that Bryant's equations determine a relative algebroid by taking the functions $c^i_{jk}, \rho^\mu_i$ to be the structure functions of the relative derivation. This is exactly the process of ``forgetting about the manifold $P$" as we did in the examples of Section \ref{sec:AlgebroidsInClassificationProblems}. In \cite{Bryant2014}, Bryant discusses many examples where an existence and classification problem for a geometric structures is described by Bryant's equations, giving rise to a variety of interesting relative algebroids.

\begin{example}[Universal/tautological relative algebroid]
Let $A\to N$ be a vector bundle and let $p_1\: \DD^1_A\to N$ be the projection. A derivation on $A$ relative to $p_1$ is a section of $p_1^*\DD^1_A$, the bundle $\DD^1_A\to N$ pulled back to itself. There is always a \textbf{tautological section} $\UD$ of this bundle given by
\[
    \UD\big\vert_{\D_x} = \D_x\in (\DD^1_A)_{p_1(\D_x)}, \quad \mbox{for $\D_x\in \DD^1_A$}. 
\] 
Interpreting elements of $\DD^1_A$ as pointwise derivations (Example \ref{ex:PointwiseDerivations}), we can explicitly describe the corresponding \textbf{tautological derivation} as
\[
    (\UD \alpha)\big\vert_{\D_x} = \D_x\alpha\in \wedge^{\bullet+1} A_x^*.
\]
The resulting triple $(A, p_1, \UD)$ is called the \textbf{universal relative algebroid} for $A$. Its name is justified by the following lemma.
\begin{lemma}[\cite{FernandesSmilde2025}, Proposition 4.3]\label{lem:UniversalPropertyUniversalAlgebroid}
    Let $(A, p, \D)$ be a relative algebroid relative to $p\:M\to N$. Then there is a map $c_\D\: M\to \DD^1_A$ such that $p = p_1\circ c_\D$ and $\D = c^*_\D \UD$. Conversely, given any map $c\: M\to \DD^1_A$ such that $p= p_1\circ c$, there is a unique relative algebroid $(A, p,\D)$ with $c_\D = c$. 
\end{lemma}
The map $c_\D$, obtained from the composition $M\xrightarrow{\D} p^*\DD^1_A\to \DD^1_A$, is called the \textbf{classifying map} of $(A, p, \D)$.

Geometrically, the universal relative algebroid $(A, p_1, \UD)$ describes the classification problem for $A$-coframes without restriction. An $A$-coframe on a manifold $P$ is a bundle map $(\theta, r)\: P\to A$ that is fiberwise an isomorphism.

If $V$ is a vector space, then according to the symbol short exact sequence (\ref{eq:SymbolSES}), $\DD^1_V\cong \Hom(\wedge^2V, V)\to \{*\}$. The \textbf{tautological bracket} dual to the tautological derivation $\UD$ on $V$ relative to $p_1\: \Hom(\wedge^2V, V)\to \{*\}$ is then given by
\[
[\cdot, \cdot]^{\Breve{}}\: \wedge^2V \to C^\infty(\Hom(\wedge^2V, V), V), \quad [v, w]^{\Breve{}}(c) = c(v, w).
\]
\end{example}

\begin{example}[Surfaces with $|\nabla K|=1$]
    The equations governing the classification problem for surfaces with $|\nabla K|=1$ as discussed in Example \ref{ex:IntroSurfacesWithConstantCurvature} define an algebroid relative to the projection $p\: S^1\times \RR\to \RR$, $(\varphi, K)\mapsto K$. 
\end{example}

Relative algebroids appear in many different contexts outside of classification problems for geometric structures.

\begin{example}[Total derivative on jet spaces]\label{ex:JetSpacesSimple}
    The total differentiation of horizontal forms is an example of a relative derivation. Suppose we are looking at functions $u(x, y)$ in two variables. The first jet space $J^1(\RR^2, \RR)$ of such functions is described by promoting $u$ and its partial derivates $u_x, u_y$ to actual coordinates. Thus, $J^1(\RR^2, \RR)$ is diffeomorphic to $\RR^5$ with coordinates $(u_x, u_y, u, x, y)$. There are projections 
    \begin{align*}
        \pi_1&\: J^1(\RR^2, \RR)\to \RR^2, & p_1&\: J^1(\RR^2, \RR)\to \RR\times \RR^2, & \pi&\:\RR\times \RR^2\to \RR^2
    \end{align*}
    The horizontal forms on $\RR\times \RR^2 = J^0(\RR^2, \RR)$ are precisely $\Omega^\bullet(\pi^*T\RR^2)$. For instance, a horizontal one-form $\alpha$ is given by
    \[
        \alpha = f(u, x, y) \d x + g(u, x, y)\d y
    \]
    for some functions $f, g\in C^\infty(J^0(\RR^2, \RR))$. 

    The total (horizontal) derivative $\D$ is the derivation on $\pi^*T\RR^2$ relative to $p_1$ determined by
    \[
        \left\{
            \begin{aligned}
                &\D (\d x) =0, &  &\D (\d y) = 0, \\ 
                &\D x = \d x, & &\D y = \d y, \\
                &\D u = u_x\d x + u_y \d y. & & 
            \end{aligned}
        \right.
    \]
    Consequently, the total derivative on $f = f(u, x, y)$ and $\alpha = f \d x +g \d y$ are given by
    \begin{align*}
        \D f &= \left(\del_x f + (\del_u f) u_x\right) \d x+ \left(\del_y f + (\del_u f)u_y\right)\d y, \\
        \D \alpha &= \left( \del_x f + (\del_u f) u_x - (\del_y f + (\del_u f)u_y\right) \d x \wedge \d y.
    \end{align*}
    Thus, horizontal differentiation defines a relative algebroid $(\pi^*T\RR^2, p_1, \D)$. We will return to this example in more generality in Section \ref{sec:RelativeAlgebroidsAndPDEs}.
\end{example}

\begin{example}[Control systems]
    Relative derivations appear everywhere in control theory. A simple control system is given by output space $\RR^n$ and a control space $\RR^k$. Given a path $y(t) = (y^\varrho(t))$ in the control space, the output $x(t) = (x^\mu(t))$ is often determined through an ordinary differential equation \cite{HermannKerner1977}:
    \[
        \dot{x}(t) = f(x(t); y(t)).
    \]
    We can rewrite this ODE as a coframe problem for intervals as
    \[
        \left\{
            \begin{aligned}
                &\D(\d t) = 0,\\
                &\D x^\mu = f^\mu(x; y) \d t.
            \end{aligned}
        \right.
    \]  
    If we forget about the underling interval, then the resulting structure consists of
    \begin{itemize}
        \item a vector bundle $\underline{\RR}\to \RR^n$ over the output space with dual basis $\theta$,
        \item a derivation $\D$ relative to the projection $\RR^k\times \RR^n\to \RR^n$ determined by
        \[
            \left\{
                \begin{aligned}
                    &\D \theta = 0, \\
                    &\D x^\mu = f^\mu(x;y) \theta
                \end{aligned}
            \right.
        \]
    \end{itemize}
    Concepts such as controllability \cite{HermannKerner1977} can be reformulated and generalized to relative algebroids.
\end{example}

\subsection{Prolongation}
The exterior derivative on a manifold always squares to zero. Therefore, for ordinary algebroids $(A, \D)$ to admit realizations, the differential $\D$ must square to zero. For relative algebroids this requirement no longer makes sense on the nose, and the consequences of the de Rham differential squaring to zero are more subtle. To understand them better, we introduce the following concept.

\begin{definition}
    Let $(A, p, \D)$ be an algebroid relative to $p\: M\to N$ with $B:=p^*A$. A \textbf{prolongation} of $(A, p, \D)$ is an algebroid $(B, p_1, \D_1)$ relative to $p_1\: M_1\to M$ such that
    \begin{itemize}
        \item $\D_1$ \textbf{extends} $\D$:
        \[
            \D_1\circ p^* = p_1^*\circ \D,
        \]
        \item $\D_1$ \textbf{completes} $\D$:
        \[
            \D_1 \circ \D = 0.
        \]
    \end{itemize}
\end{definition}
We depict a prolongation by the following diagram, where $B_1 = p_1^*B$.
\[
    \xymatrix{
        B_1 \ar[r] \ar[d] & B \ar[r] \ar[d] \ar@/^1pc/@{-->}[l]^{\D_1} & A \ar[d] \ar@/^1pc/@{-->}[l]^{\D} \\
        M_1 \ar[r]^{p_1} & M \ar[r]^{p_1} & N
    }
\]

\subsubsection{Computing prolongations}\label{subsubsec:ComputingProlongations}

Suppose a relative algebroid $(A, p, \D)$ in a local (dual) frame $(\theta^i)$ and submersion coordinates $(y^\varrho, x^\mu)$ is described by structure functions $c^{i}_{jk}, \rho^\mu_i$. Then an \textit{extension} $\D_1$ of $\D$ is determined by how it acts on the fiber coordinates $y^\varrho$. Thus we can add $\D_1$ to the structure equations to obtain
\begin{equation}\label{eq:ExtentionInCoordinates}
    \left\{
    \begin{aligned}
        &\D \theta^i = -\frac{1}{2} c^i_{jk} \theta^i\wedge\theta^k, \\
        & \D x^\mu = \rho^\mu_i\theta^i,\\
        & {\color{magenta}\D_1 y^\varrho = \xi^\varrho_i \theta^i}.
    \end{aligned}
    \right.
\end{equation}
We can treat the $\xi^\varrho_i$'s as new variables (as opposed to functions on $M$) and impose the completion condition $\D_1\circ \D=0$ to obtain equations for $\xi^\varrho_i$. This carves out a space $M_1$ on which $\D_1$ defines a derivation relative to the projection $M_1\to M$. The object $\D_1\circ \D$ can be interpreted as the \emph{torsion} of the relative algebroid (see \cite{FernandesSmilde2025}).

\begin{example}
    For this example we let $A = \underline{\RR^2}\to \RR^2 =\{x, y\}$, the trivial bundle with fiber $\RR^2$ over the space $\RR^2$. Let $\D$ be a derivation relative to to the projection $\RR\times \RR^2 \to \RR^2$, $(z, x, y)\mapsto (x, y)$ determined by
    \[
        \left\{
            \begin{aligned}
                &\D \theta^1 = \theta^1\wedge\theta^2, \\
                &\D \theta^2 = 0, \\
                &\D x= z \theta^1, \\
                &\D y = z \theta^2.
            \end{aligned}
        \right.
    \]
    An extension $\D_1$ is determined by
    \[
        {\color{magenta}\D_1 z = z_1 \theta^1 + z_2 \theta^2}
    \]  
    for some new variables $z_1, z_2$. Imposing $\D_1\circ \D =0$ gives four equations: $\D_1\circ \D \theta^i =0$ (which are automatic) and 
    \begin{align*}
        \D_1\circ \D x &= (\D_1z)\wedge \theta^1 + z\theta^1\wedge\theta^2 = 0,\\
        \D_1\circ \D y &= (\D_1 z) \wedge\theta^2 =0.
    \end{align*}
    It follows that $z_1=0$ and $z_2=z$, so $\D_1$ is determined by
    \[
        {\color{magenta}\D_1 z = z\theta^2}.
    \]Thus, a prolongation of $(A,p, \D)$ is the algebroid relative to the identity $(\underline{\RR^2},\D_1)\to \RR\times \RR^2$. In this case, one can check that $\D_1^2=0$ so $(\underline{\RR^2}, \D_1)$ is an actual \emph{Lie algebroid}, isomorphic to the action algebroid $\RR^2\ltimes \RR^3$ where $\RR^2$ is equipped with the non-abelian Lie algebra structure.
\end{example}

\begin{example}\label{ex:NoProlongation}
    Even if there is a prolongation, the prolongation itself might not admit another prolongation. Let us modify the example above, so that we get a derivation $\D$ on $\underline{\RR^2}\to \RR^2$ relative to the projection $(z, x, y)\mapsto (x, y)$ defined by
    \[
        \left\{
            \begin{aligned}
                &\D \theta^1 = \theta^1\wedge\theta^2, \\
                &\D \theta^2 = 0, \\
                &\D x= z \theta^1, \\
                &\D y = z \theta^2 {\color{green}\ -\ \theta^1}.
            \end{aligned}
        \right.
    \]
    Following the same procedure as in the previous example, we find that the prolongation is still a derivation on $\underline{\RR^2}\to \RR^3$ relative to the identity but now determined by
    \[
        {\color{magenta}\D_1 = z \theta^2} {\color{green} \ + \ \theta^1}.
    \]
    However, there is no prolongation of the prolongation $(\underline{\RR^2}, \D_1)$, because
    \[
        \D_1^2z = 2\theta^1\wedge\theta^2 \neq 0. 
    \]
    The object $\D_1^2 z$ is related to the \emph{intrinsic curvature} of the relative algebroid (see \cite{FernandesSmilde2025}).
\end{example}

\begin{example}
    In the examples above, the equations on the extra variables $\xi^\varrho_i$ came purely from the symbol of the derivation $\D$. This doesn't always have to be the case. Here is a simple example where the equations are derived from the $c^i_{jk}$'s. 

    Consider $\underline{\RR^3}\to \{*\}$ and the derivation relative to $\RR=\{x\} \to \{*\}$ defined by
    \[
        \left\{
            \begin{aligned}
                &\D \theta^1 = \theta^2\wedge \theta^3 + x\theta^1\wedge\theta^2, \\
                &\D \theta^2 = \theta^3 \wedge\theta^1 + x \theta^2 \wedge \theta^3, \\
                &\D \theta^3 = \theta^1 \wedge\theta^2 + x \theta^3 \wedge \theta^1.
            \end{aligned}
        \right.
    \]
    Then an extension $\D_1$ is determined by
    \[
        {\color{magenta}\D_1 x = x_1 \theta^1 + x_2 \theta^2 + x_3\theta^3}.
    \]
    The $\D_1\circ \D = 0$ implies
    \begin{align*}
        \D_1 \circ \D \theta^1 = (\D_1 x)\wedge \theta^1\wedge\theta^2 = x_3\theta^3\wedge\theta^1\wedge\theta^2
    \end{align*}
    and similar for $\theta^2, \theta^3$. This implies that the only extension must satisfy $\D_1 x = 0$. This time, $\D^2_1x=0$, and the resulting Lie algebroid is a bundle of Lie algebras.
\end{example}

\begin{example}
    The above examples are relative algebroids \textit{of finite type}. This is typically not the case. Here is an example of an infinite type relative algebroid. Start with a derivation on $\underline{\RR^2}\to \RR\cong \{x_0\}$ relative to the projection $(x_1, x_0)\mapsto x_0$ determined by
    \[
        \left\{
            \begin{aligned}
                &\D \theta^1 = 0,\\
                &\D \theta^2 = - f(x_0)\theta^1\wedge\theta^2,\\
                &\D x_0 = x_1\theta^2
            \end{aligned}
        \right.
    \]
    for some function $f$ of one variable. The prolongation $\D_1$ is now determined by
    \[
        \D_1 \circ \D x_0 = (\D_1x_1)\wedge \theta^2 - x_1f(x_0)\theta^1\wedge\theta^2=0.
    \]
    It follows that
    \[
        {\color{magenta} \D_1 x_1 = x_1f(x_0)\theta^1 + x_2\theta^2}.
    \]
    For the next prolongation, we must have
    \[
        \D_2\circ \D_1 x_1 = \D_2(x_1f(x_0)\theta^1 + x_2\theta^2) = -(2x_2 f(x_0) + x_1^2f'(x_0))\theta^1\wedge\theta^1 + (\D_2 x_2)\wedge\theta^2=0,
    \]  
    so
    \[
        {\color{magenta} \D_2 x_2 = (2x_2f(x_0) + x_1^2f'(x_0))\theta^1 + x_3\theta^2}.
    \]
    All the higher prolongations $\D_k x_k$ are of the form
    \[
        \D_kx_k = g_k(x_0, \dots, x_k) \theta^1 + x_{k+1}\theta^2,
    \]
    where $g_k(x_0, \dots, x_k)$ is function. There is a recursive formula for $g_k$ obtained from
    \begin{align*}
        \D_{k+1}\circ \D_k x_k &= D_{k+1}\left( g_k \theta^1 + x_{k+1}\theta^2\right),\\
        &= -\left( \sum_{i=0}^k \dfrac{\del g_k}{\del x_i} x_{i+1} + x_{k+1}f(x_0)\right)\theta^1\wedge\theta^2 + (\D_{k+1}x_{k+1})\wedge\theta^2.
    \end{align*}
    It follows that
    \begin{align*}
        g_{k+1}(x_0, \dots, x_{k+1}) &= \sum_{i=0}^k \dfrac{\del g_k}{\del x_i} x_{i+1} + x_{k+1}f(x_0),
    \end{align*}
    so in particular, the function $g_k$ is polynomial in the variables $x_1, \dots, x_k$. If $f(x_0)$ it itself a polynomial, then $g_k$ is a polynomial in all its variables.
\end{example}

\subsubsection{The prolongation of a relative algebroid}

The algorithm to for building the prolongation of $(A, p, \D)$ as in the examples above relies on two steps. The first is to define a \textit{space of extensions}, and the second is to carve the \textit{space of completions} out of the space of extensions. On the space of completions we will define a canonical relative algebroid that is a prolongation of $(A, p, \D)$.

Let $(A, p, \D)$ be a relative algebroid. Recall from Section \ref{subsec:StructureOfRelativeDerivations} that there is a map $p_*\: \DD^1_B\to p^*\DD^1_A$ with kernel $\Hom(B, \ker Tp)$. The \textbf{space of (pointwise) extensions} is the affine subbundle $L\subseteq \DD^1_B$ modeled on $\Hom(B, \ker Tp)$ defined by
\[
        L := p_*^{-1}(\im \D) = \left\{ \tilde{\D}_x\in \DD^1_B \ | \ p_*(\tilde{\D}_x) = \D_x \mbox{ for all $x\in M$}\right\}.
\]
The space of extensions comes naturally with an inclusion into $i_L\: L \hookrightarrow \DD^1_B$ intertwining the bundle projections $p_1$. An element $\tilde{\D}_x\in L$ is determined by a (pointwise) lift of the anchor
\[
    \xymatrix{
        {} & T_x M \ar[d]\\
        B_x \ar[ur]^-{\tilde{\rho}_x} \ar[r]_-{\rho_x} & T_{p(x)}N
    }
\]
which is equivalent to defining $\tilde{\D}_x$ on the fiber coordinates as in Equation (\ref{eq:ExtentionInCoordinates}). The lift of the anchor $\tilde{\rho}_x$ is the symbol of a derivation $\tilde{\D}_x\in L$.

By Lemma \ref{lem:UniversalPropertyUniversalAlgebroid}, the universal derivation $\UD$ on $B$ relative to $p_1$ restricts to a relative derivation $\tilde{\D}_1$ on $L$ relative to $p_1$! Explicitly, 
\begin{equation}
    {\color{magenta}{(\D_1\alpha)\big\vert_{\tilde{D}_x} = \tilde{\D}_x \alpha}}, \quad \mbox{for $\alpha\in \Omega^1(B)$ and $\tilde{\D}_x\in L$.}
\end{equation}
This is the coordinate-free equivalent of the derivation defined by Equation (\ref{eq:ExtentionInCoordinates}), with the parameters $\xi^\varrho_i$ treated as coordinates. The resulting relative algebroid $(B, p_1, \tilde{\D}_1)$ is an extension of $(A, p, \D)$ by construction, but not necessarily a completion.

To obtain an actual prolongation, we have to restrict to the \textbf{space of completions}, defined by
\[
    M^{(1)} = \left\{ \tilde{\D}_x\in L \ | \ \tilde{\D}_x \circ \D = 0 \right\} 
\]
where we conveniently interpret elements in $\tilde{\D}_x\in L \subseteq \DD^1_B$ as pointwise derivations (Example \ref{ex:PointwiseDerivations}) to make sense of the composition $\tilde{\D}_x\circ \D$. The space $M^{(1)}$ may not be smooth (in fact, it may be empty), as Example \ref{ex:NoProlongation} shows. 

\begin{definition}
    A relative algebroid $(A,p, \D)$ is called \textbf{1-integrable} if $M^{(1)}$ is smooth and the projection $p_1\: M^{(1)}\to M$ is a submersion.

    In that case, the \textbf{prolongation} of $(A, p, \D)$ is the relative algebroid $(A, p, \D)^{(1)}:=(B, p_1, \D^{(1)})$ with as classifying map $c_{\D^{(1)}}\: M^{(1)}\hookrightarrow \DD^1_B$ the inclusion.
\end{definition}

The prolongation $(B^{(1)}, p_1, \D^{(1)})$ is universal among all prolongations in the following sense, even though it's not always smooth.

\begin{lemma}
    Let $(A, p,\D)$ be a relative algebroid. An algebroid $(B, p_1, \D_1)$ relative to $p_1\: M_1\to M$ is a prolongation if and only if its the classifying map $c_{\D_1}\: M_1\to \DD^1_B$ has image in $M^{(1)}$. 
\end{lemma}

\subsection{Formal integrability}

Higher prolongations of a relative algebroid can be defined iteratively (if they are smooth). We call a relative algebroid $(A, p, \D)$ \textbf{$k$-integrable} when it is $(k-1)$-integrable and its $(k-1)$-th prolongation is 1-integrable. In that case, the $k$-th prolongation is the prolongation of the $(k-1)$-th prolongation:
\[
    (B^{(k-1)}, p_k, \D^{(k)}) = (B^{(k-2)}, p_{k-1}, \D^{(k-1)})^{(1)}.
\]
\begin{definition}
    A relative algebroid is \textbf{formally integrable} when it is $k$-integrable for all $k$. A formally integrable relative algebroid is called a \textbf{relative Lie algebroid}.
\end{definition}
For a relative Lie algebroid, the prolongations assemble into a tower
\[
    \xymatrix{
    {\left(B^{(\infty)}, \D^{(\infty)}\right)\:\ldots} \ar[r] \ar@<-12pt>[d] & B^{(k)} \ar[d] \ar[r] &  \ldots \ar[r] \ar@{-->}@/^1pc/[l]^{\D^{(k)}} & B^{(1)} \ar[d] \ar[r] & B \ar@{-->}@/^1pc/[l]^{\D^{(1)}} \ar[r] \ar[d] & A \ar[d] \ar@{-->}@/^1pc/[l]^{\D} \\
    M^{(\infty)}\:\ldots \ar[r]                                            & M^{(k)} \ar[r]_{p_k}    & \ldots \ar[r]                                                               & M^{(1)} \ar[r]_{p_1} & M \ar[r]_{p}    & N.
    }               
\]
The resulting derivation $\D^{(\infty)}$ squares to zero on \emph{profinite} sections of the profinite vector bundle $B^{(\infty)}\to M^{(\infty)}$, so $(B^{(\infty)}, \D^{(\infty)})$ is a \emph{profinite Lie algebroid}, by which we mean an algebroid object in the category of profinite manifolds (as opposed to a pro-object in the category of Lie algebroids). Note that while $M^{(\infty)}$ is typically infinite-dimensional, the vector bundle $B^{(\infty)}\to M^{(\infty)}$ has finite rank.

\begin{remark} 
    Checking whether each prolongation is 1-integrable is often intractable. In practice, due to a version of Goldschmidt's formal integrability criterion for relative algebroids \cite[Theorem 3.10]{FernandesSmilde2025}, it is possible to check for formal integrability just by checking for 1-integrability and involutivity of a \emph{tableau of derivations} associated to the relative algebroid. See \cite{FernandesSmilde2025} for more details.
\end{remark}

\subsubsection{Existence of realizations}

Robert Bryant proved an existence result for Bryant's Equations \ref{eq:BryantsEquations} using the Cartan-K\"{a}hler theorem for exterior differential systems \cite[Theorem 3 and 4]{Bryant2014}. It translates in our language to an existence result for realizations of relative Lie algebroids.

\begin{theorem}
    Let $(A, p, \D)$ be a \emph{real analytic} relative Lie algebroid. Then for each $k$ and $x\in M^{(k)}$, there exists a realization $(P, \theta_k, r_k)$ of $(B^{(k-1)}, p_k, \D^{(k)})$ such that $x\in \im r_k$. 
\end{theorem}

\begin{question}[Extension problem]
    A different way of approaching the realization problem for relative algebroids is to start with a relative algebroid $(A, p, \D)$ and look for an (local) extension $\tilde{D}\in \Gamma(\DD^1_B)$ of $\D$ such that $\tilde{\D}^2=0$. We call this the \textbf{extension problem} for relative algebroids. 
    
    While a solution to the extension problem does not give any further information about the pointwise existence, it does give (locally) an actual Lie algebroid $(B, \tilde{\D})$. A (local) groupoid integrating this Lie algebroid gives rise to large families of solutions, and the integration theory of Lie algebroids gives obstructions to the existence of complete solutions, or rather, local solutions that are complete in the isotropy direction.

    We expect that the formal theory of the extension problem is equivalent to that of the relative algebroid, and that formal integrability together with real analyticity implies existence of local solutions to
\end{question}

\subsection{An example revisited}

Let us now make the theory explicit with the example of Riemannian surfaces with $|\nabla K|=1$ (Example \ref{ex:IntroSurfacesWith|nablaK|=1}). We refer to \cite[Section 6]{FernandesSmilde2025} for more details. The equations for the classification problem for such structures gives rise to an algebroid $(\underline{\RR^3}, p, \D)$ relative to $p\:S^1\times \RR\to\RR$, $(\varphi, K)\mapsto K$, with derivation determined by
\[
    \left\{
    \begin{aligned}
        &\D \theta^1 = - \theta^3 \wedge \theta^2, \\
        &\D \theta^2 = \theta^3 \wedge \theta^1, \\
        &\D \theta^3 = K \theta^1\wedge\theta^2,\\
        &\D K = (\cos \varphi)\theta^1 + (\sin \varphi)\theta^2. 
    \end{aligned}
    \right.
\]
As described in Section \ref{subsubsec:ComputingProlongations}, the prolongation can be computed by first setting 
\[
    \D^{(1)}\varphi = \varphi_1\theta^1 + \varphi_2\theta^2 + \varphi_3\theta^3,
\]
and then subjecting $\varphi_i$ to the equations $\D_1\circ \D K =0$. In this case, it follows that
\[
    \D^{(1)}\varphi = - \theta^3 + c_1 (\del_\varphi \D K)
\]
where $c_1$ parametrizes the fiber coordinates of the first prolongation space and $\del_\varphi \D K$ is the shorthand notation for 
\[
    \del_\varphi\D K = -(\sin \varphi)\theta^1 + (\cos\varphi)\theta^2.
\]
This defines the first prolongation as the algebroid $(\underline{\RR^{(3)}}, p_1, \D^{(1)})$ relative to $p_1\: \RR\times S^1\times \RR\to S^1\times \RR$, $(c_1, \varphi, K) \mapsto (\varphi, K)$. 

Higher prolongations can be explicitly computed, and are of the form 
\begin{equation}\label{eq:RevisitedProlongation}
    \D^{(k+1)}c_{k} = f_k[c_1, \dots, c_k] \D K + c_{k+1}(\del_\varphi \D K),
\end{equation}
for some explicit polynomial $f_k[c_1 \dots, c_k]$. Thus, the relative algebroid is formally integrable, and determines a profinite relative algebroid 
\[
    (\underline{\RR^3}, \D^{(\infty)}) \to \RR^{\infty}\times S^1\times \RR,
\]
where $\RR^{\infty}\times S^1\times \RR$ has coordinates $(\ldots, c_3, c_2, c_1, \varphi, K)$.

\subsubsection{Analysis of the relative algebroid}
    There are several interesting ways in which the profinite algebroid leads to insights into the classification problem of these structures.
    
\subsubsection*{Solutions with additional symmetry}    The rank of the anchor of $(\underline{\RR^3}, \D^{(\infty)})$ is generically 3, but drops over the locus $\{ \ldots = c_3 = c_2 = 0\}$. This is an invariant submanifold for $(\underline{\RR^3}, \D^{(\infty)}$, so it restricts to a finite-dimensional Lie algebroid $(\underline{\RR^3}, \D)$ over $\RR\times S^1\times \RR$, with derivation defined by
    \[
        \left\{
        \begin{aligned}
            &\D \theta^1 = - \theta^3\wedge\theta^2, & & & &\D K = (\cos \varphi) \theta^1 + (\sin \varphi) \theta^2\\
            &\D \theta^2 = \theta^3\wedge\theta^1, & & & &\D \varphi = \theta^3 + c_1 (\del_\varphi \D K)\\
            &\D \theta^3 = K\theta^1\wedge\theta^2, & & & &\D c_1 = -(c_1^2 + K) \D K.
        \end{aligned}
        \right.
    \]
    This describes the realization problem of surfaces with $|\nabla K| = 1$ and with one-dimensional isometry group. The Lie algebroid is integrable, and an integration can be explicitly given in terms of a solution to a Riccati equation. All non-extendable, simply connected solutions are contractible, and have isometry group $\mathbb{Z}_2 \ltimes \RR$. Other non-extendable realizations can be obtained by choosing a lattice in $\RR$ and taking the quotient, resulting in solutions diffeomorphic to cylinders with isometry group $\mathbb{Z}_2\ltimes S^1$. 

    \subsubsection*{The profinite Lie algebroid}$\ $ Equation (\ref{eq:RevisitedProlongation}) suggests that the algebroid structure is better described by taking a different frame. If $(e_1, e_2, e_3)$ is the basis of dual to $\theta^1, \theta^2, \theta^3$, we could take
    a new basis that decouples the relative part of the profinite
    \begin{equation}
        \begin{aligned}
            b_1 &= (\cos \varphi) e_1 + (\sin \varphi) e_2,\\
            b_2 &= -(\sin \varphi) e_1 + (\cos \varphi) e_2 - c_1 e_3, \\
            b_3 &= e_3
        \end{aligned}
    \end{equation}
    The only non-zero bracket between these sections is
    \begin{align*}
        [b_1, b_2] &= - c_0 b_2,
    \end{align*}
    while the anchors are given by
    \begin{align*}
        \rho(b_1) &=  \del_K + \sum_{k=1}^\infty f_k \del_{c_k}, & \rho(b_2) &= \sum_{k=1}^\infty c_{k+1} \del_{c_k}, & \rho(b_3) &= \del_\varphi. 
    \end{align*}
    Remarkably, only one of those vector fields is of true profinite nature, the other ones are (levelwise) honest vector fields. Integral curves of $\rho(b_2)$ are determined by the choice of a function of one variable, while the integral curves of $\rho(b_1)$ and $\rho(b_3)$ are determined by their initial point. Hence, using flows it is possible to construct realizations of this profinite algebroid by starting with an integral curve of $\rho(b_2$ and flowing along $\rho(b_1)$ and $\rho(b_3)$

    This case is rather special because there is only one vector field that is truly of profinite nature. This is because the only non-zero Cartan character of the system is $s_1=1$. We expect that realizations of relative algebroids with $s_1=1$ can always be constructed in a similar manner. This would generalize a Lie theorem for PDEs of class $\omega=1$ to relative algebroids \cite{Kruglikov2012}.
    
    \subsubsection*{Analytic integrations} 
    The profinite algebroid $(\underline{\RR}^3, \D^{(\infty)}) \to \RR^\infty \times S^1\times \RR$ has odd behaviour. The leaves of the equivalence relation generated by points connected through realizations are infinite-dimensional, for instance, even though the realizations themselves have finite dimension. One way to get closer to more intuitive behavior of the relative algebroid is to change the resolution to the prolongation tower and restrict to \emph{real analytic} realizations only. The real analytic realization correspond precisely to the the resolution of convergent power series $\RR^\omega \times S^1\times \RR$, where
    \[
        \RR^\omega = \{ (\ldots, c_3, c_2, c_1) \ | \ (c_k)_k \mbox{ defines a convergent power series} \}.
    \]
    The relative algebroid restricts to $\RR^\omega \times S^1\times \RR$, where it has unique local analytic realization passing through each point. In particular, it has well-defined leaves The algebroid can be restricted to each leaf where it becomes a finite-dimensional Lie algebroid. The global solutions to this realization problem can then be computed explicitly using techniques from \cite{FernandesStruchiner2019}. Computing the leaves explicitly is rather difficult, as they depend on analytic continuation as well as a maximal solution to a Riccati ODE. We expect that there could be a smooth groupoid integrating $(\underline{\RR^3}, \D^{(\infty)})\to \RR^\omega\times S^1\times \RR$, but the nature of its smooth structure has not yet been investigated.

\section{PDEs and relative algebroids}\label{sec:RelativeAlgebroidsAndPDEs}

As hinted to in Example \ref{ex:JetSpacesSimple}, the total derivative on jet spaces can be interpreted as a relative derivation. Before we explain how a general PDE gives rise to a relative algebroid, we recall some facts about jet spaces and in particular the Cartan distribution.

Let $q\: Q\to X$ be a submersion. The space of $k$-jets $J^kq$ consists locally of order $k$ Taylor expansions of sections of $q$. We denote a $k$-jet of a local section $\sigma$ at $x\in X$ by $j^k_x \sigma$. There are projections
\begin{align*}
    q_k&\: J^k q\to X, & p_k&\: J^kq\to J^{k-1}q,
\end{align*}
where the first projects a jet to the point that is the origin of the Taylor expansion, and the second truncates an order $k$ Taylor expansion to order $k-1$.

Given a local section $\sigma$ of $q$, its $k$-jet $j^k\sigma$ defines a section of $J^k q$. Such sections of $q_k$ are called \textbf{holonomic}.

First jets of sections can alternatively be described by a complement $H$ to $\ker Tq\subset TQ$, or as splittings of $Tq$. We will frequently use the following identification:
\begin{equation}\label{eq:FirstJetSpaceSplittings}
    J^1q\cong \left\{ h_x\: T_{q(x)}X \to T_x Q\ | \ Tq\circ h_x = \id_{T_{q(x)}X} \right\} \subset \Hom(q^*TX, TQ). 
\end{equation}
Another useful fact is that the space of $k$-jets $J^kq$ naturally sits inside the first jet space of $q_{k-1}$:
\[
    J^kq \hookrightarrow J^1q_{k-1}, \quad j^k_x\sigma \mapsto j^1_x (j^{k-1}\sigma).
\]

The (relative) Cartan distribution, the central structural object of a jet space, detects when sections of $q_k$ are holonomic. On the level of first jets it can be described as the subbundle
\[
    \Cc \subset p_1^*TQ, \quad \Cc_{h_x} = \im h_x
\]
for $h_x\in J^1q$. It is actually a complement to $p_1^*\ker Tq$ and therefore the map $p_1^*Tq\: p_1^*TQ\to p_1^*q^*TX$ restricts to an isomorphism $\Cc\cong p_1^*q^*TX$. The Cartan distribution on $J^k q$ can be obtained by restricting the Cartan distribution on $J^1q_{k-1}$ along the inclusion $J^kq\hookrightarrow J^1q_{k-1}$. 

A section $\tau$ of $q_{k}$ is holonomic if and only if it is tangent to the Cartan distribution:
\[
    \im T_{\tau(x)}p_1\circ T_x\tau \subseteq \Cc_{\tau(x)},
\]
for all $x$ in the domain of $\tau$.

The space of first jets can also be described using the language of derivations! To explain how, recall first that the de Rham differential $\d$ is a section of $\DD^1_{TX}$. Also, there is a natural map $q_*\: \DD^1_{q^*TX}\to q^* \DD^1_{TX}$ that fits together with the symbol map into the following diagram (Section \ref{subsec:StructureOfRelativeDerivations}):
\[
    \xymatrix{
        \DD^1_{q^*TX} \ar[r]^-\sigma \ar[d]_{q_*} & \Hom(q^*TX, TQ) \ar[d]^{q_*} \\
        q^*\DD^1_{TX} \ar[r]^-\sigma & q_1^*\Hom(TX, TX).
    }
\]
The symbol of the de Rham differential $\d$ is $\id_{TX}\in \Hom(TX, TX)$, and it is clear from Equation (\ref{eq:FirstJetSpaceSplittings}) that $J^1q = q_*^{-1}(\id_{TX})$, so it follows that
\begin{equation}
    J^1q\cong q_*^{-1}(\im q^*\d) = \{ \D_x\in \DD^1_{q^*TX} \ | \ q_*(\D_x) = \d_{q(x)}\} \subset \DD^1_{q^*TX}.
\end{equation}

\subsection{PDEs and prolongation}

Let $q\: Q\to X$ be a submersion. A \textbf{partial differential equation (PDE)} of order $k$ on $q$ is a subspace $E_k\subset J^kq$. A \textbf{solution} to $E_k$ is a section $\sigma$ of $q$ such that $\im j^k\sigma\subset E_k$. The terminology for prolongation and formal integrability actually comes from the theory of formal PDEs developed by Spencer and his school \cite{Goldschmidt1967, Spencer1969}, so it makes sense to very briefly recall it in order to make the connection with relative algebroid more clear later. A great modern reference is the book by Seiler \cite{Seiler2010}.

The \textbf{first prolongation} of a PDE $E_k \subset J^kq$ is the PDE of order $k+1$ defined by
\[
    E^{(1)}_k := J^1(q_k\big\vert_{E_k}) \cap J^{k+1}q \subset J^1 q_k.
\]
If it is smooth, and $p_{k+1}\: E^{(1)}_k\to E_k$ is a surjective submersion, we say that $E_k$ is \textbf{1-integrable}. It is $l$-integrable when it is $(l-1)$-integrable $E^{(l-1)}_k$ is 1-integrable and we define the $l$-th prolongation by $E^{(l)}_k = (E^{(l-1)}_k)^{(1)}$. A PDE is \textbf{formally} integrable if it is $l$-integrable for all $l$. 

Since the $l$-th prolongation of $E_k$ is a PDE of order $k+l$, we will also use the following notation:
\[
    E_{k+l} = E^{(l)}_k\subset J^{k+l}q.
\]  

\subsection{PDEs with symmetries}

One of the biggest advantages of relative algebroids over PDEs is that relative algebroids are stable under quotients. Here, we only consider point symmetries, although the statements are valid in the more general context of Pfaffian fibrations and their symmetries \cite{Smilde2025}.

Symmetries of PDEs are often described by pseudogroups. A \textbf{pseudogroup} for us means a collection of local diffeomorphisms $\Gamma\subseteq \Diffloc(X)$ on a manifold $X$ that is closed under composition, inversion, restriction and gluing. 

If $\Gamma$ is a pseudogroup on $P$, then its collection of $k$-jets
\[
 J^k\Gamma:= \{ j^k_x\varphi \ | \ \varphi\in \Gamma, x\in \dom \varphi\}.
\]  
is a groupoid on $P$, with composition $j^k_{\varphi(x)}\psi \cdot j^k_x \varphi = j^k_x(\psi\circ\varphi)$. For the local diffeomorphisms, we have
\[
    J^k\Diffloc(X) = J^k(X\times X), 
\]
where $J^k(X\times X)$ is the groupoid of $k$-jets of bisections of the pair groupoid $X\times X\toto X$. Note that if $\Gamma$ is a pseudogroup on $X$, then $J^k\Gamma\subseteq J^k\Diffloc(X)$, and this can be interpreted as a PDE of order $k$ on (bi)sections of $X\times X$. 

\begin{example}
    Let $q\: Q\to X$ be a submersion. We denote by $\Gamma_q$ the pseudogroup of diffeomorphisms preserving $\ker Tq$. Locally, an element $\varphi\in \Gamma_q$ covers a local diffeomorphism $\phi$ on $X$, i.e. $q\circ \varphi = \phi \circ q$ locally. Because of this property, germs of elements of $\Gamma_q$ act on germs of sections of $q$ through
    \[
        \varphi_* \sigma = \varphi \circ \sigma \circ \phi^{-1}.
    \]
    This action descends to the level of $k$-jets, so every $\varphi\in \Gamma_q$ lifts to a local diffeomorphism $\varphi^{(k)} \: J^kq\dashrightarrow J^kq$, called the $k$-th prolongation of $\varphi$. 
\end{example}

If $E_k\subset J^kq$ is a PDE, the \textbf{point symmetries} of $E_k$ are those $\varphi\in \Gamma_q$ such that $\varphi^{(k)}(E_k)\subseteq E_k$. 

In practice, one is either interested in computing the full group of point symmetries of a PDE, or the PDE itself already comes naturally with a (sub)set of point symmetries. In the latter case, if $\Gamma\subset \Gamma_q$ is a pseudogroup, we say that a PDE $E_k\subset J^kE$ is $\Gamma$-invariant when each element of $\Gamma$ is a point symmetry of $E_k$, i.e. $\varphi^{(k)}(E_k)\subseteq E_k$ for all $\varphi\in \Gamma$.

Predetermined symmetries of a PDE often appear together with natural bundles. If $\Gamma$ is a pseudogroup on $X$, a \textbf{$\Gamma$-natural bundle} is a submersion $q\: Q\to X$ with a lift $L\:\Gamma\to \Gamma_q$ that respects inversion and composition and such that each lift $L(\varphi)$ covers $\varphi$ on $X$ \cite{Accornero2021, PalaisTerng1977}. For example, metrics on $X$ are sections of $q\:S^2T^*X\to X$, which is a $\Diffloc(X)$-natural bundle. So, a PDE $E_k\subset J^kq$ on metrics on $X$ is diffeomorphism-invariant when it is invariant under (the lift of) the action of $\Diffloc(X)$. It is possible that the PDE posses more symmetries.

\subsection{The relative algebroid of a PDE}\label{sec:TheRelativeAlgebroidOfAPDE}

Next, we describe how a PDE naturally gives rise to a relative algebroid. Let $q\: Q\to X$ be a submersion and $E_k\subset J^kq$ a PDE. We assume for simplicity that $E_{k-1}:=p_k(E_k)\subset J^{k-1}q$ is smooth and that $p_k$ restricts to a submersion $p_k\: E_k\to E_{k-1}$. This is not strictly necessary, but keeps us in the realm of algebroids relative to submersions, rather then foliations. For a treatment without this assumption, we refer to \cite{FernandesSmilde2025, Smilde2025}. We will also assume that $E_k\subset E^{(1)}_{k-1}$. This assumption implies that the submersion $p_k\: E_k\to E_{k-1}$ factors through the inclusion $E_k\hookrightarrow J^1q_{k-1}$, where $q_{k-1}\: E_{k-1}\to X$ is the projection.

The relative algebroid $(q^*_{k-1}TX, p_k, \D^{E_k})$ underlying $E_k$ is defined through the classifying map $c\: E_k \to \DD^1_{q_{k-1}^*TX}$ given by $E_k \hookrightarrow J^1q_{k-1} \hookrightarrow \DD^1_{q_{k-1}^*TX}$. 

More concretely, the derivation is determined by
\begin{align}
    \D^{E_k} q_{k-1}^*\alpha &= q_k^*\d \alpha, & \langle \D^{E_k} f\big\vert_{h_x}, v\rangle&= \langle \d f, h_x(v)\rangle,
\end{align}
for $\alpha\in \Omega^\bullet_{TX}$, $f\in C^\infty(E_{k-1})$, $h_x\in E_k\subset J^1q_{k-1}$ and $v\in T_{q_{k-1}(x)}X$.

This construction completely embeds the formal theory of PDEs into that of relative algebroids.

\begin{theorem}[\cite{FernandesSmilde2025, Smilde2025}]
    Let $E_k\subset J^kq$ be a PDE with underlying relative algebroid $(q_{k-1}^*TX, p_k, \D^{E_k})$. Then:
    \begin{enumerate}
        \item The prolongation space of $E_k$ is canonically isomorphic to the prolongation space of $(q_{k-1}^*TX, p_k, \D^{E_k})$. In particular, the PDE $E_k$ is 1-integrable of and only if $(q_{k-1}^*TX, p_k, \D^{E_k})$ is.
        \item If $E_k$ is 1-integrable, the prolongation of $(q_{k-1}^*TX, p_k, \D^{E_k})$ is canonically isomorphic to the relative algebroid underlying $E^{(1)}_k$. 
        \item Germs of solutions to $E_k$ are in 1-1 correspondence with realizations of $(q_{k-1}^*TX, p_k, \D^{E_k})$ modulo germs of diffeomorphisms.
    \end{enumerate}
\end{theorem}

\subsection{Relative algebroids and PDEs with symmetries}

The relative algebroid is completely natural with respect to point symmetries of a PDE. 

\begin{lemma}[\cite{Smilde2025}]
    Let $\varphi$ be a point symmetry of a PDE $E_k\subseteq J^kq$. Then $\varphi$ induces an isomorphism on the underlying relative algebroid $(q_{k-1}^*TX, p_k, \D^{E_k})$, i.e. a bundle isomorphism
    \[
        \varphi_*\: q_k^*TX\to q_k^*TX
    \]
    covering $\varphi_*\: q_{k-1}^*TX\to q_{k-1}^*TX$ such that $\varphi^*\circ \D^{E_k} = \D^{E_k}\circ \varphi^*$.       
\end{lemma}
The map $\varphi_*$ is defind as follows. Since $\varphi$ preserves $\ker Tq$, it induces a bundle isomorphism $\varphi_*\: q^*TX\to q^*TX$. Together with $\varphi^{(k)}\: E_k\to E_k$ it induces an isomorphism $\varphi_*$ of $q_k^*TX$.

The groupoid language provides a convenient way to deal with quotients by pseudogroups. Let $\Gamma\subset \Gamma_q$ be a pseudogroup of point symmetries of a PDE $E_k\subset J^kq$. Since $\varphi^{(k)}(x)$ depends only on the $k$-jet of $\varphi$, there is a groupoid action
\[
\xymatrix{
	J^k\Gamma \ar@<0.5ex>[d] \ar@<-0.5ex>[d] & E_k \ar@`{[l]+/l-0.2pc/+/d+1.8pc/,[l]+/l-0.2pc/+/u+1.8pc/}[0, 0]_{m} \ar[ld]^{\mu_k}    \\
	Q & 
}
\]
where $\mu\: J^kq\to Q$ sends the $k$-jet $j^k_x\sigma$ to its target $\sigma(x)$. The space $E_{k-1}= p_k(E_k)$ is also $\Gamma$-invariant, and the projection
\[
    (J^k\Gamma \curvearrowright E_k) \xrightarrow{p_k} (J^{k-1}\Gamma \curvearrowright E_{k-1})
\]
is equivariant with respect to $J^{k}\Gamma\to J^{k-1}\Gamma$. 

For $k\geq 1$, there is also a natural representation
\[
    J^k\Gamma\curvearrowright q_k^*TX
\]
covering the action of $J^k\Gamma$ on $E_k$ (so technically it is a representation of the action groupoid $J^k\Gamma\ltimes E_k$).

For the action on of $J^k\Gamma$ on the bundle of derivations, we have to distinguish between two cases.

\textbf{Case $k>1$}$\ $ In this case there is also a representation $J^{k-1}\Gamma\curvearrowright q_{k-1}^*TX$. Together with the $J^k\Gamma$ action on $E_k$ it induces a representation on the bundle of relative derivations
\begin{equation}\label{eq:RepresentationOnDerivations}
    J^k\Gamma \curvearrowright p_k^*\DD^1_{q_{k-1}^*TX}, \quad j^k_{\mu_k(x)}\varphi \cdot \D_x = \varphi_*\circ \D_x \circ \varphi^*, 
\end{equation}
where we regard pointwise derivation $\D_x\in p_k^*\DD^1_{q_{k-1}^*TX}$ as maps (Example \ref{ex:PointwiseDerivations})
\[
     \D_x\: \Omega^\bullet\left(q_{k-1}^*TX\right)\to \wedge^{\bullet+1}T_{p_k(x)}^*X.
\]
Note that Equation (\ref{eq:RepresentationOnDerivations}) depends exactly on the $k$-jet of $\varphi$, because $\varphi^*\: \Omega^\bullet_{q_{k-1}^*TX} \to \Omega^\bullet_{q_{k-1}^*TX}$ depends on the $(k-1)$-jet of $\varphi$.

\textbf{Case $k=1$}$\ $ In this case, there is no representation $J^0\Gamma$ on $q^*TX$. However, Equation (\ref{eq:RepresentationOnDerivations}) is still valid for $k=2$, so in this case there is a representation
\[
    J^{k+1}\Gamma\curvearrowright p^*_k \DD^1_{q^*_{k-1}TX}.
\]
In other words, in both cases there is a representation of $J^{k+1}\Gamma$ on the bundle of relative derivations, but only when $k>1$, it descends to a represenation of $J^k\Gamma$. 

\begin{theorem}[\cite{Smilde2025}]
    Let $E_k\subseteq J^kq$ be a PDE and $\Gamma\subseteq \Gamma_q$ a pseudogroup of point symmetries of $E_k$. Then the derivation $D^{E_k}$ of its underlying relative algebroid is an invariant section of for the representation $J^{k+1}\Gamma\curvearrowright p_k^*\DD^1_{q_{k-1}^*TX}$. 
\end{theorem}

\begin{corollary}\label{cor:InvariantSectionsPreserved}
    If $E_k\subseteq J^kq$ is a PDE with symmetry $\Gamma$, then the derivation $\D^{E_k}$ preserves $\Gamma$-invariant sections:
    \[
        \D^{E_k}\: \Omega^\bullet\left(q_{k-1}^*TX\right)^\Gamma \to \Omega^\bullet\left(q_{k-1}^*TX\right)^\Gamma 
    \]
\end{corollary}
As a consequence, $J^{k-1}\Gamma$ and $J^k\Gamma$ are Lie groupoids, and the actions $J^{k-1}\Gamma\curvearrowright E_{k-1}$ and $J^{k}\Gamma\curvearrowright E_k$ are free and proper, then the derivation descends to a derivation on $q_{k-1}^*TX/\Gamma$ relative to $E_k/\Gamma\to E_{k-1}/\Gamma$. There is a lot of potential for applications of this result, particularly in the relation with work by Olver and his collaborators on pseudogroups and invariants \cite{OlverPohjanpelto2005, OlverPohjanpelto2009}. For example, in \cite{olverPohjanpelto2012} it is shown that if $J^k\Gamma\curvearrowright E_k$ is free, then so is the action $J^{k+1}\Gamma\curvearrowright E^{(1)}_k$.

\subsection{Application: the relative algebroid of a relative distribution}
We end the notes with an application of taking a quotient of a PDE by symmetries to get a relative algebroid.

In the literature of PDEs, it is often customary to consider equations on arbitrary submanifolds of a fixed dimension, rather than sections of a submersion (see e.g. \cite{PuglieseSparanoVitagliano2023}). In this setup, we consider the space of $k$-jets of $n$-dimensional submanifolds of a fixed ambient manifold $Q$, and denote it by $J^k(Q, n)$. Note that $J^1(Q, n)$ is the Grassmanian of $n$-planes in $TQ$, so it comes with a tautological bundle $\tau_n\to J^1(Q, n)$ that embeds canonically into $p_1^*TQ$. 

If $Q$ happens to be the total space of a submersion $q\: Q\to X$ with $n$-dimensional base, then $J^kq$ is an open and dense subset of $J^k(Q, n)$. In fact, this can always be done locally, so $J^k(Q, n)$ can be covered by opens of the form $J^kq$ where $q\: Q\dashrightarrow X$ is some local submersion. The Cartan distribution on $J^1q$ corresponds to the tautological bundle $\tau_n\to J^1(Q, n)$ restricted to $J^1q$. 

We have seen in Section \ref{sec:TheRelativeAlgebroidOfAPDE} that there is a canonical structure of a relative algebroid on $q_k^*TX\to J^kq$ for any submersion $q\: Q\to X$. Its construction relied on the de Rham differential on $X$ (or the involutive complement $\ker Tq$ to the Cartan distribution \cite{FernandesSmilde2025, Smilde2025}), which does not generalize automatically to $J^k(Q, n)$. Instead, it is possible to use a quotient construction to obtain a relative algebroid over $J^k(Q, n)$ whose realizations correspond to submanifolds of $Q$. 

\newcommand{\imm}{\mathrm{imm}}
To explain how, we consider $J^k_{\imm} (\RR^n, Q)$ space of $k$-th order jets of (local) immersions $\RR^n\to Q$ (we could also take any other $n$-dimensional manifold $X$ instead of $\RR^n$). Note that $J^1_\imm(\RR^n, Q)$ is a PDE on the sections of the projection $s\: \RR^n\times Q\to \RR^n$, and that $J^k_\imm(\RR^n, Q)$ coincides with the $(k-1)$-th prolongation of $J^1_\imm(\RR^n, Q)$. 

The natural $\Diffloc(\RR^n)$-action on $\RR^n\times Q$ preserves the fibers of the projection, and induces point symmetries of $J^1_\imm(\RR^n, Q)$, giving rise to natural actions $J^k\Diffloc(\RR^n)\curvearrowright J^k_\imm(\RR^n, Q)$. 

The algebroid relative to $p_k\: J^k_\imm(\RR^n, Q)\to J^{k-1}_\imm(\RR^n, Q)$ has as base vector bundle $s_{k-1}^*T\RR^n$ with a natural $J^{k-1}\Diffloc(\RR^n)$ action as long as $k-1\geq 1$. We have the following quotients:
\begin{align*}
    J^k_\imm(\RR^n, Q)/J^k\Diffloc(\RR^n)&\cong J^k(Q, n), & s^*_k T\RR^n/J^k\Diffloc(\RR^n) &\cong p_{k, 1}^* \tau_n,
\end{align*}
for $k\geq 1$, where $p_{k,1}\: J^k(Q, n)\to J^1(Q, n)$ is the projection. 

By Corollary \ref{cor:InvariantSectionsPreserved}, the derivation $\D^{k+1}$ relative to $p_k\: J^k_\imm (\RR^n, Q)\to J^{k-1}_\imm(\RR^n, Q)$ preserves $\Diffloc(\RR^n)$-invariant sections and therefore descends to a derivation $\D^k$ on $p_{k-1,1}^*\tau_n$ relative to $p_k\: J^k(Q, n)\to J^{k-1}(Q, n)$.

\begin{proposition}\label{prop:RelativeAlgebroidJetsOfSubmanifolds}
    For $k\geq 1$, there is a natural structure of a relative algebroid on $p_{k-1, 1}^*\tau_n$ relative to $p_{k+1}\: J^{k+1}(Q, n)\to J^k(Q, n)$ with derivation $\D^{k+1}$. It has the following properties:
    \begin{enumerate}
        \item The derivation $\D^{k+1}$ is a prolongation of $\D^k$.
        \item The quotient maps $s_k^*T\RR^n\to p^*_{k,1}\tau_n$ are morphisms of relative algebroids.
        \item Realizations of $\D^k$ correspond to immersions of rank $n$ into $Q$.
    \end{enumerate}
\end{proposition}

\begin{remark}
    Even when considering the more general setting of algebroids relative to foliations, the bundle $\tau_n$ itself does not seem to carry a structure of a relative algebroid for which the realizations are the immersions. The quotient does not even exist in the category of bundles with partial flat connections. While the quotient map $J^1_\imm(\RR^n, Q)\to J^1(Q, p)$ is a map of foliations, there is no flat $\ker Tp_1$-connection on $\tau_n$ induced through the quotient. The reason is that the only $\Diffloc(\RR^n)$-invariant \emph{flat} connection of $(s_1^*T\RR^n, \onabla)\to J^1_\imm(\RR^n,Q)$ is the zero section.
\end{remark}

It is also possible to write down an explicit formula for the relative derivation of the underlying relative algebroid. 

\begin{lemma}\label{lem:FormulaDerivationJetsOfSubmanifolds}
    The derivation of the relative algebroid 
    \[
        \xymatrix{
            p_2^*\tau_n \ar[r] \ar[d] & \tau_n \ar[d] \ar@/^1pc/@{-->}[l]^{\D} \\ 
            J^2(Q, n) \ar[r]^{p_2} & J^1(Q, n)
        }
    \]
    is given by
    \[
        \left(\D \alpha\right)_{j^2_xi}(v_0, \dots, v_k) = \left( \d (j^1i)^*\alpha\right)\left( \left(T_xi\right)^{-1}(v_0), \dots, \left( T_xi\right)^{-1}(v_k)\right)
    \]
    for $\alpha\in \Omega^\bullet(\tau_n)$ and $v_0, \ldots, v_n\in (\tau_n)_{\im T_xi}$. 

    The bracket dual to this derivation is determined by
    \[
        [X, Y]_{j^2_xi} = T_x i\left( [i^*X, i^*Y]_x\right),
    \]
    for $X, Y\in \Gamma(\tau_n)$.
\end{lemma} 
Note that if $i\:N\to Q$ is an immersion, then $((j^1i)^*\alpha)_x = (T_xi)^*\alpha_{j^1_x i}$ is a well-defined honest differential form on $N$ because $(\tau_n)_{j^1_xi} = \im T_x i$. It depends pointwise on the first jet of $i$, so the expression $\d (j^1i)^*\alpha$ makes sense and depends pointwise on the second jet of $i$. Similarly, the expression $(i^*X)_{\tilde{x}} = \left(T_{\tilde{x}}i\right)^{-1}(X_{\im T_{\tilde{x}} i})$ defines an honest vector field on $N$.

\begin{proof}[Sketch of the proof of Lemma \ref{lem:FormulaDerivationJetsOfSubmanifolds}]
    There are two ingredients for this proof. The first ingredient is that quotient map $\pi\: s^*T\RR^n\to p_2^*\tau_n$ is a morphism of relative algebroids, meaning that $\pi^*\circ \D = \D \circ \pi^*$. The second is more general and stems from the fact that the algebroid relative to $J^2_\imm(\RR^n, Q)\to J^1_\imm(\RR^n, Q)$ is the prolongation of the canonical algebroid relative to $J^1_{\imm}(\RR^n, Q)\to \RR^n\times Q$. Generally, if $(A, p, \D)$ is a relative algebroid and 
    \[
        \xymatrix{
        B^{(1)} \ar[r] \ar[d] & B \ar[r] \ar[d] \ar@/^1pc/@{-->}[l]^{\D^{(1)}} & A \ar[d] \ar@/^1pc/@{-->}[l]^{\D} \\
        M^{(1)} \ar[r]^{p_1} & M \ar[r]^{p_1} & N
    } 
    \]
    is the prolongation, then any realization $(P, \theta, r)$ differentiates to a realization $(P, \theta^{(1)}, r^{(1)})$ where $r^{(1)}\: P\to M^{(1)}\subset \DD^1_B$ is given by
    \[
        r^{(1)}(x) = (\theta_x^{-1})^*\circ \d_x \circ \theta^*.
    \]
    The tautological nature of $\D^{(1)}$ forces it to be given pointwise by 
    \[
        \D^{(1)}_{r^{(1)}(x)} = (\theta_x^{-1})^*\circ \d_x \circ \theta^*.
    \]

   Since every immersion is a realization for the algebroid relative to $J^1_{\imm}(\RR^n, Q) \to \RR^n\times Q$, we can apply this fact to the prolongation sequence $J^2_\imm(\RR^n, Q)\to J^1_\imm (\RR^n, Q)\to \RR^n\times Q$ and use hat the quotient map is a morphism of relative algebroids to derive the desired formula.
\end{proof}

\subsubsection{Relative distributions} Let $Q$ be any manifold and $p\: M\to Q$ a map. A \textbf{distribution on $Q$ relative to $p$} is a subbundle
\[
    H\subseteq p^*TQ.
\]
A distribution of rank $n$ relative to $p$ is characterized by the map $h\: M\to J^1(Q, n)$ with the property that $H_{x} = (\tau_n^Q)_{h(x)}$ for all $x\in M$. This property implies that there is a canonical identification $H\cong h^*\tau_n^Q$. We use $\tau^Q_n\to J^1(Q, n)$ for the tautological bundle over the Grassmanians of $n$-planes in the manifold $Q$. 

An \textbf{integral manifold} of a distribution relative $H\subset p^*TQ$ relative to $p\: M\to Q$ a submanifold $i\: N\hookrightarrow Q$ with a lift $\tilde{i}\:N\to M$ such that $\im T_xi = H_{\tilde{i}(x)}$. If $p$ is a submersion, then this forces $\tilde{i}$ to be a submanifold as well. 

The \textbf{partial prolongation} of a distribution $H\subset p^*TQ$ relative to $p\: M\to Q$ is the space
\[
    J^1_HM := \left\{ \tilde{H}_x\in J^1(M, n)\ | \ T_xp(\tilde{H}_x) = H_x\subset T_{p(x)} Q \right\}.
\]
The inclusion $\tilde{h}\: J^1_HM\hookrightarrow J^1(M, n)$ is now the classifying map of a distribution $\tilde{H}\subset p_1^*TM$ relative to $p_1\: J^1_HM\to M$.

A relative distribution is enough to define a relative anchor but not a relative bracket. This requires information of the second jet. To extract this information, we use the classifying map $h\: M\to J^1(Q, n)$. This differentiates to a map
\[
    h_*\: J^1_HM\to J^1(J^1(Q, n), n), \quad h_*(\tilde{H}_x) = T_xh (\tilde{H}_x)\subset T_{h(x)} J^1(Q, n).
\]
It is well-defined because $Tp_1\left(h_*(\tilde{H}_x)\right) = H_x\subset T_{p(x)} Q$ is an $n$-dimensional subspace.

The \textbf{first prolongation} of the distribution $H$ relative to $p\: M\to N$ is the space
\[
    M^{(1)} = h_*^{-1}\left( J^2(Q, n)\right)
\]
If smooth, there is an algebroid $h^*\tau_n^Q\to M$ relative to $p_1\: M^{(1)}\to M$ whose derivation is determined by the symbol coming from the canonical inclusion
\[
    p_1^*h^*\tau_n^Q\cong \tilde{h}^*\tau^M_n \hookrightarrow p_1^*TM,
\]
and the requirement that
\[
    \xymatrix{
        p_1^*h^*\tau_n^Q \ar[r] \ar[d] & p_2^*\tau^Q_n \ar[d]\\
        M^{(1)} \ar[r]^{h_*} & J^2(Q, n)
    }
\]
is a morphism of relative algebroids.

\bibliographystyle{abbrv}
\bibliography{bib}

\end{document}